\documentclass[12pt]{amsart}
\usepackage[centertags]{amsmath}

\usepackage{amsfonts}

\usepackage{amssymb}

\title[Automorphism Groups of  Models of PA]{Automorphism Groups of 
Countable \\ Arithmetically Saturated Models of \\ Peano Arithmetic}
\author{James H. Schmerl}
\date{\today}
\def\into{\longrightarrow}

\def\harp{\hspace{-4pt} \upharpoonright \hspace{-4pt}}
\def\pa{{\sf PA}}
\def\pas{\pa^*}
\def\ta{{\sf TA}}
\def\wt{\widetilde}

\newcommand{\MM}{{\mathcal M}}
\newcommand{\NN}{{\mathcal N}}

 \DeclareMathOperator{\scl}{Scl}
 \DeclareMathOperator{\rep}{Rep}
  \DeclareMathOperator{\ssy}{SSy}  
 
  \DeclareMathOperator{\aut}{Aut}
 
  \DeclareMathOperator{\Th}{Th} 
  \DeclareMathOperator{\tp}{tp}
  \DeclareMathOperator{\Def}{Def}
   \DeclareMathOperator{\igap}{igap}
    
   \DeclareMathOperator{\eq}{Eq}

\def\ls{{\mathcal L}\raisebox{.2ex}{*}}

  \DeclareMathOperator{\Lt}{Lt}

\begin{document}

\begin{abstract} If $\MM, \NN$ are countable, arithmetically saturated models of Peano Arithmetic and $\aut(\MM) \cong \aut(\NN)$, then the  Turing-jumps of 
$\Th(\MM)$ and $\Th(\NN)$ are  recursively equivalent. \end{abstract}

\maketitle


\bigskip


Since 1991, when the question

\begin{quote} 
{\em Are there countable, recursively saturated models $\MM$, $\NN$ of \pa\ 
such that $\aut(\MM) \not\cong \aut(\NN)$ $($as abstract groups})?
\end{quote}
appeared in \cite{kkk}, it has been of interest to determine to what extent  
(the isomorphism type of) the group $\aut(\MM)$ of all automorphisms of 
 a countable, recursively saturated model $\MM$ of Peano Arithmetic  determines (the isomorphism type  of)  
$\MM$. It was proved in \cite{kkk}  that 
whenever  both $\MM$ and $\NN$ are countable, recursively saturated 
models of \pa\ and exactly one of them is arithmetically saturated, then 
$\aut(\MM)$ and $\aut(\NN)$ are not isomorphic {\em as topological groups}. 
 In 1994, Lascar \cite{la} 
proved  that countable, arithmetically saturated models of \pa\ have the small index property, and that result then implied that $\aut(\MM) \not\cong \aut(\NN)$ {\em as abstract groups}. 
This gave the first positive answer to the above question. A neater way, in which the use of the small index property is masked,  that  automorphism 
groups distinguish those  models that are arithmetically saturated from the other countable 
recursively saturated ones was obtained the next year in \cite[Coro.\@ 3.9]{ks2} (or see  \cite[Th.\@ 9.3.10]{ksbook}): 
{\em If $\MM$ is countable and recursively saturated, then $\MM$ is arithmetically saturated 
iff the cofinality of $\aut(\MM)$ is uncountable.} Finally, we mention that Kaye's Theorem 
\cite{kaye2} (see \S1) 
characterizing the closed normal subgroups of $\aut(\MM)$, which appeared in the same 
volume \cite{auto} as did Lascar's Theorem, yields that whenever 
$\MM, \NN$ are countable, arithmetically saturated models and $\MM$ is a model of True Arithmetic $(\ta)$ while $\NN$ is not, then $\aut(\MM) \not \cong \aut(\NN)$. 

Recall that a countable, recursively saturated model $\MM$ of $\pa$ is determined up to isomorphism by two invariants: 
its   standard system $\ssy(\MM)$ and its first-order theory $\Th(\MM)$. 
Correspondingly, there are the following complementary questions for a countable, recursively saturated 
model $\MM$ of $\pa$.

\begin{center} 

 {\em To what extent does 
$\aut(\MM)$ determine \ $\ssy(\MM)$}? \ 
 {\em $\Th(\MM)$}?

\end{center} 
Subsequent to  Lascar's  proof about the small index property, the focus has been almost entirely on   countable, arithmetically saturated models. The ``$\ssy$'' 
question  for these models was answered soon thereafter.

\bigskip 

{\sc Theorem 1}: (Kossak-Schmerl \cite{ks}) {\em If $\MM, \NN$ are countable, arithmetically saturated models of $\pa$ 
such that  $\aut(\MM) \cong \aut(\NN)$, then $\ssy(\MM) =  \ssy (\NN)$.} 

\bigskip

Progress on the ``$\Th$'' question has been slower. Nurkhaidarov \cite{n}  proved the following in 2006.

\bigskip

{\sc Theorem 2}: (Nurkhaidarov \cite{n}) {\em There are countable, arithmetically saturated models 
$\MM_0$, $\MM_1$, $\MM_2$, $\MM_3$ of $\pa$ such that whenever 
$ i < j < 4$, then $\ssy(\MM_i) = \ssy(\MM_j)$ and $\aut(\MM_i) \not\cong \aut(\MM_j)$.}

\bigskip

Although not explicitly stated in \cite{n}, the proof of Theorem~2  also proves the following stronger 
result.

\bigskip

{\sc Theorem 3}: (Nurkhaidarov \cite{n}) {\em There are completions 
$T_0$, $T_1$, $T_2$, $T_3$ of $\pa$ such that whenever 
$ i < j <4$ and $\MM_i$, $\MM_j$ are countable, arithmetically saturated models 
of $T_i, T_j$, respectively,  then $\aut(\MM_i) \not\cong \aut(\MM_j)$.}

\bigskip

Theorem~3 implies Theorem~2 because whenever ${\mathcal T}$ is a countable set of 
completions of \pa, then there is an ${\mathfrak X}$ such that each $T \in {\mathcal T}$ 
has a countable, arithmetically saturated model whose standard system is~${\mathfrak X}$.  

This paper  improves Theorem~3 by increasing the cardinal number~$4$ in that theorem  to the maximum possible  of $2^{\aleph_0}$. If $X,Y \subseteq \omega$, 
then we write $X \leq_T Y$ if $X$ is Turing-reducible (or recursive relative) to $Y$, and $X \equiv_T Y$ if $X$ is recursively equivalent to $T$ (that is, $X \leq_T Y \leq_T X$).  
As usual, the  Turing-jump of $X$ is $X'$.   The following  theorem is our principal new result.

\bigskip

{\sc Theorem 4}:   {\em If $\MM, \NN$ are countable, arithmetically saturated 
models of $\pa$ 
and  $\aut(\MM) \cong \aut(\NN)$, then $ \Th(\MM)' \equiv_T \Th(\NN)'$.} 

\bigskip

A consequence of this theorem is that the cardinal 4 in Theorem~2 can be increased to 
the maximum possible of $\aleph_0$. In fact, we get the following corollary that yields   some answers  to Question~15 in \cite[Chap.~12]{ksbook}. 

\bigskip

{\sc Corollary 5}: {\em For any   countable jump ideal ${\mathfrak X}$, 
there are  infinitely many countable arithmetically saturated 
models $\MM_0, \MM_1, \MM_2, \ldots$ of \ \pa\ such that whenever $i < j < \omega$, then $\aut(\MM_i) \not\cong \aut(\MM_j)$ and $\ssy(\MM_i) = {\mathfrak X}$.}

\bigskip

One may wonder whether Theorems~1 and 4 tell the whole story. In other words, 
if $\MM, \NN$ are countable, arithmetically saturated models of \pa\ such that $\aut(\MM) \cong \aut(\NN)$, then is it the case that $\ssy(\MM) = \ssy(\NN)$ and $\Th(\MM)' \equiv_T \Th(\NN)'$? We easily see that this is not so since 
   the 4 theories in Theorem~3 can be chosen to be recursively equivalent.
However, we can do even better.

\bigskip

{\sc Theorem 6}:  { \em For each $n < \omega$, there are  recursively equivalent completions $T_0,T_1, \ldots, T_n$ of $\pa$ such that whenever 
$i < j \leq n$ and $\MM_i,\MM_j$ are countable, arithmetically saturated models of 
$T_i,T_j$, respectively, then  $\aut(\MM_i) \not\cong \aut(\MM_j)$.}

\bigskip

The results presented in this paper suggest the question that is dual to the one asked in 
\cite{kkk} and could have just as easily been asked there.

\bigskip

{\it Question} 7: Are there countable, recursively saturated models $\MM$, $\NN$ of \pa\ 
such that $\MM \not\cong \NN$ and $\aut(\MM) \cong \aut(\NN)$? 

\bigskip

The analogous question for countable, arithmetically saturated models 
is also open.

There are 6 sections that follow this introductory one. Some preliminaries are in 
\S1, which consists of some notation, definitions and results that will be used in 
the succeeding sections. The next two sections do not overtly refer to automorphism groups.
Some special types of models of \pa, the lofty models and those with the $\omega$-property, are discussed in \S2. In \S3, we consider substructure lattices, carefully reviewing 
a result from \cite{ks2}. Our main result, Theorem~4, is proved in \S4. Theorem~6 is proved in \S5, and some additional results are given in \S6. 

\smallskip

Roman Kossak and Ermek Nurkhaidarov are thanked for their helpful 
comments on various precursors of this paper.

\bigskip


{\bf  \S1.\@ Some Preliminaries.} Notation and terminology used here will generally follow 
\cite{ksbook}. The reader should refer to \cite{ksbook} for insufficiently explained notions. 
The proof of Theorem~4 relies on a number of results 
that are proved or stated in \cite{ns}. It is suggested that the reader have that paper available.

The language appropriate for \pa\ is ${\mathcal L}_\pa = \{+,\times, \leq 0,1\}$.
It is to be tacitly understood that all models  referred to in this paper are models of $\pa$.
All models are assumed to have the standard model ${\mathbb N} = (\omega,+,\times, \leq ,0,1)$ as a submodel. Models will be denoted by (possibly adorned) script letters 
such as $\MM, \NN, \MM_1, \ldots$, and their universes are denoted by the corresponding roman letters
$M,N,M_1, \ldots$,  although models and their universes  may  occasionally be confounded. 

Suppose that $X,Y \subseteq \omega$. As already  mentioned, 
$X \leq_T Y$ iff $X$ is Turing-reducible to $Y$, $X \equiv Y$ iff $X$ and $Y$ are  recursively equivalent,  and $X'$ is the Turing-jump 
of $X$. For a small ordinal $\alpha$, 
$X^{(\alpha)}$ is the $\alpha$-th jump of $X$. 
 If there is $n < \omega$ such that $X \leq_T Y^{(n)}$, then $X$ is {\it arithmetically reducible} to $Y$ (in symbols: $X \leq_a Y$).    If $X \leq_a Y \leq_a X$, then $X$ and $Y$ are {\it arithmetically equivalent} 
(in symbols: $X \equiv_a Y$). We let $X\oplus Y = \{2x : x \in X\} \cup \{2y+1 : y \in Y\}$.

A {\it Turing ideal} is a nonempty subset ${\mathfrak X} \subseteq {\mathcal  P}(\omega)$ such that whenever $X,Y \in {\mathfrak X}$ and $Z \leq_T X \oplus Y$, then $Z \in {\mathfrak X}$. A {\it jump ideal} is a Turing ideal ${\mathfrak X}$ such that $Y \in {\mathfrak X}$ whenever $Y \leq_a X \in {\mathfrak X}$. 

Suppose that $\MM$ is an arbitrary model. If $A \subseteq M$, then the model 
generated by $A$, denoted by $\scl(A)$, is the smallest elementary substructure 
of $\MM$ containing $A$.   The model $\MM$ is finitely generated iff 
$\MM = \scl(a) \ (= \scl(\{a\}))$ for some $a \in M$. A subset $I$  is a cut of $\MM$ iff $0 \in I$ and whenever $x \leq y \in I$, then $x +1 \in I$. A cut $I$ is invariant iff $I = \sup(I \cap \scl(0))$ or $I = \inf(\scl(0) \backslash I)$, and 
it is exponentially closed iff $2^x \in I$ whenever $x \in I$.  We let $\Lt(\MM)$ be the lattice of elementary substructures of $\MM$ and $\Lt_0(\MM)$ be its $\vee$-subsemilattice consisting of those models in $\Lt(\MM)$ that are 
finitely generated. It is a consequence of Ehrenfeucht's Lemma (\cite[Theorem~1.7.2]{ksbook}) that 
$\aut(\MM)$ is trivial whenever $\MM$ is finitely generated.
The standard system of $\MM$ 
is $\ssy(\MM)= \{D \cap \omega : D$ is a definable subset of $\MM\}$. In general,  $\ssy(\MM)$ is a Scott set or, equivalently, 
$({\mathbb N}, \ssy(\MM)) \models {\sf WKL}_0$. 
  If $\MM$ is recursively saturated, 
then $\MM$ is arithmetically saturated iff $({\mathbb N}, \ssy(\MM)) \models {\sf ACA}_0$ 
iff $({\mathbb N}, \ssy(\MM)) \models {\sf RT}^3_2$ iff $\ssy(\MM)$ is a jump ideal 
iff $\omega$ is a strong cut of $\MM$. Here, we are letting ${\sf RT}^n_2$ denote infinite Ramsey's Theorem for 2-colored $n$-sets. 

The usual interval notation will be used. If $\MM$ is a model and $a,b \in M$, 
then $[a,b] = \{x \in M : a \leq x \leq b\}$ and  $[a,b) = \{x \in M : a \leq x < b\}$.

If $T \supseteq \pa$ is a theory (which, for us, is  a (not necessarily deductively closed) consistent set of sentences) and $X \subseteq \omega$, then $X$ is a  real represented by $T$ if there is a unary formula $\varphi(x)$ such that for each $n < \omega$,
$$
n \in X \Longleftrightarrow T \vdash \varphi(n) \Longleftrightarrow T \not\vdash \neg\varphi(n).
$$
As usual, $\rep(T)$ is the set of reals represented by $T$. If $T$ is complete,
then $\rep(T)$  
is the standard system of the prime model of $T$.

If $G =  \aut(\MM)$ and $A \subseteq M$, then the pointwise stabilizer of $A$ is the subgroup
$G_{(A)} = \{g \in G : g(a) = a$ for all $a \in A\}$ and the setwise stabilizer is 
$G_{\{A\}} = \{g \in G : g[A] = A\}$. If $A = \{a\}$, then $G_a = G_{(A)} = G_{\{A\}}$. When considering $G$ as a topological group,  the stabilizers of finite subsets of $M$ are its basic open subgroups. Equivalently, the basic open subgroups are the 
pointwise stabilizers of finitely generated elementary submodels. Since finitely generated models 
do not have any nontrivial automorphisms, the basic open subgroups are also 
the 
setwise stabilizers of finitely generated elementary submodels.

The following theorem has already been mentioned.

\bigskip

{\sc Kaye's Theorem:} {\em If $\MM$ is a countable recursively saturated model and 
$H \leq G = \aut(\MM)$, then the following are equivalent$:$

\begin{itemize}

\item[(1)] $H$ is a closed normal subgroup of $G$.

\item[(2)] $H = G_{(I)}$, where $I \subseteq M$ is an invariant, exponentially closed cut.

\end{itemize}}

\bigskip

Even though the next theorem will be not explicitly used in this paper, we state it since it 
shows that for the implication $(1) \Longrightarrow (2)$ in Kaye's Theorem there is a unique such $I$. 

\bigskip

{\sc Smory\'{n}ski's Theorem}: {\em Suppose that $\MM$ is a countable recursively saturated model and that 
$H \leq G = \aut(\MM)$. If  $H = G_{(I)}$ for some cut $I \subseteq M$, then there is a unique exponentially closed cut $J \subseteq M$ such that $H = G_{(J)}$.}

\bigskip

Smory\'{n}ski actually proved more (see \cite[Theorem~8.4.2]{ksbook}): 
If $\MM$ is a countable recursively saturated model and $I \subseteq M$ is an exponentially closed cut, then there is $f \in \aut(\MM)$ such that 
$I = I_{\it fix}(f) = \{x \in M : f(y) = y$ for all $y < x\}$. 
 
 For the record, we state Lascar's Theorem on the small index property.
 
 \bigskip
 
 {\sc Lascar's Theorem}: {\em Suppose that $\MM$ is a countable, arithmetically saturated model. If $H \leq \aut(\MM)$, then $H$ is open iff its index 
 $|\aut(\MM) : H|$ is countable.}
 
 \bigskip

Scott \cite{sco} introduced the notion of a Scott set and proved two related theorems 
concerning  standard systems and sets of represented reals. We will need  the following variants of 
these two theorems, the first of which also appears as \cite[Theorem~13.6]{kaye}.

\bigskip 

{\sc Theorem 1.1}: (Wilmers \cite{wil}) {\em Suppose that $T$ is a completion 
of \pa\ and ${\mathfrak X}$ is a  set of subsets of $\omega$. The following are equivalent\hspace{1pt}$:$ 

$(1)$ ${\mathfrak X}$ is a countable Scott set and $T \in {\mathfrak X}$.

$(2)$ There is a countable recursively saturated $\MM \models T$ such that $\ssy(\MM) = {\mathfrak X}$.}

\bigskip

 If $X \subseteq \omega$ and $i < \omega$, then $(X)_i = \{j < \omega : \pmb{\langle} i,j \pmb{\rangle} \in X\}$. If ${\mathfrak X}$ is a set of subsets of $\omega$ and $X \subseteq \omega$, then ${\mathfrak X}$ is {\it enumerated} by $X$ 
 if ${\mathfrak X} = \{(X)_i : i < \omega\}$. 
 
 \bigskip

{\sc Theorem 1.2}: (Knight \cite[Coro.\@ 1.6]{knight} and D.\@ Marker) {\em If ${\mathfrak X}$ is a Scott set enumerated by $X$ and $X \leq_T Y$, then $\pa$ has  a completion  $T$ such that  $\rep(T) = {\mathfrak X}$ 
 and $T \equiv_T Y$.}
 
 \bigskip
 
 The next proposition collects together some well known properties of countable, recursively saturated 
 models most of which can be found in various places in \cite{ksbook}.

\bigskip

{\sc Proposition 1.3}: {\em Suppose that $\MM$ is a countable, recursively saturated model. 

\begin{itemize}

\item[(1)] $\MM$ is tall.

\item[(2)] $\MM$ is generated by a set of indiscernibles of any countable order type with no last element.

\item[(3)] $\MM$ is homogeneous.

\item[(4)] If $\MM_0 \succ_{\sf cf} \MM $, then $\MM_0$ is recursively saturated.
\item[(5)] There is $\MM_0 \prec_{\sf end} \MM$ such that $\MM_0 \cong \MM$.

\end{itemize}}

\bigskip

What follows in this paragraph is a small exception to our convention that all models 
considered here are models of \pa. Consider the language  $\ls = {\mathcal L}_\pa \cup \{U\}$,  where $U$ is 
 a new unary relation symbol, and then let $\pas$ be the ${\mathcal  L}^*$-theory obtained from \pa\ by adjoining all instances of the induction scheme in this expanded language. Each model of $\pas$ expands a model of \pa. Every statement 
in this paper that applies to models of \pa\ has a natural extension that applies to  
models of $\pas$. We will have several occasions when we will want to refer to 
such an extended version of some result, and we will do so by referring to its $*$-version.

\bigskip


{\bf \S2.\@ Loftiness and the $\omega$-property.}  This section is concerned with  some 
properties of models that were introduced in \cite{ks4}, \cite{ks5} and \cite{om1}. 
Results of this section will be used in the proofs of the main results although 
automorphism groups do not appear here. The results may 
 be of independent interest.   Theorem~2.8 has to do with constructing 
countable models that have the $\omega$-property but are not recursively saturated. Corollary~2.9 is a common generalization of Proposition~1.3(2),(4) and (5). We begin with the definitions.

\bigskip

{\sc Definition 2.1:} Suppose that $\MM$ is a nonstandard model and $I$ is a cut.

\smallskip

(1)  $I$  is {\em upward monotonically $\omega$-lofty} if there is $a \in M$ such that $I = \sup\{(a)_i : i < \omega\}$.

(2)  $I$  is {\em downward monotonically $\omega$-lofty} if there is $a \in M$ such that $I = \inf\{(a)_i : i < \omega\}$.

(3)  $I$  is {\em uniformly $\omega$-lofty} if there is $a \in M$ 
such that whenever $p < I < q$, then $p < (a)_i < q$ for some $ i < \omega$.

(4) $\MM$ is {\it uniformly $\omega$-lofty} if  for any $b \in M$ there is $a \in M$ such that whenever $\omega < e \in M$, then 
$\scl(b) \subseteq \{(a)_i : i < e\}$.

(5) $\MM$ has the $\omega$-{\em property} if there is $\NN \succ_{\sf end} \MM$ such that 
$M$ is an upward monotonically $\omega$-lofty cut of $\NN$. 

\bigskip

Various notions of loftiness were introduced and studied in \cite{ks4} and \cite{ks5}. 
Definitions (1), (2)  and (3) are from \cite[Def.\@ 3.1]{ks5}. Definition (4) is not the 
one given in \cite[Def.\@ 1.4(iii)]{ks4}, although it is equivalent. One direction of this equivalence is given in  \cite[Theorem~1.7(1b)]{ks2}; the other  is easy to see. 
It is straightforward to see that every recursively saturated model is uniformly $\omega$-lofty. In fact, if we let $t_0(x), t_1(x), t_2(x), \ldots$ be a recursive list of all Skolem terms, then $\MM$ is recursively saturated iff for every $b \in M$ there is $a \in M$ 
such that $\MM \models (a)_i = t_i(b)$ for all $i < \omega$ (\cite[Prop.\@ 1.6]{ks5}).
 Clearly, every uniformly $\omega$-lofty model is tall. 
Finally, (5) was introduced by Kossak in \cite{om1} and studied by him in \cite{om1} and \cite{om2}.

According to \cite{om1}, the following lemma is implicit in \cite{ks5};  it is explicitly proved 
in \cite[Theorem~2.7]{om1}. 

\bigskip

{\sc Lemma 2.2:} {\em If $\MM$ is tall and  has the $\omega$-property, then $\MM$ is 
uniformly $\omega$-lofty.}

\bigskip

{\sc Lemma 2.3:} {\em Suppose that $\MM$ is  countable and uniformly $\omega$-lofty. The following are equivalent$:$

\begin{itemize}

\item[(1)] $\MM$ is recursively saturated.

\item[(2)] $\MM$ is generated by a set of indiscernibles.

\item[(3)] $\MM$ is   isomorphic to some $\MM_0 \prec_{\sf end} \MM$.

\end{itemize}}

\bigskip

{\it Proof}.   $(1) \Longrightarrow (2)$ and $(1) \Longrightarrow (3)$ are (2) and (5) of Proposition~1.3.

For the converses, we rely on \cite[Lemma~1.8 -- Theorem~1.13]{ks5} from which it 
follows that if $\MM_0 \prec_{\sf end} \MM$ and $\MM_0$ is not recursively saturated, 
then there is $a \in M \backslash M_0$ that realizes a type not realized in $\MM_0$. 
Both of the implications $(2) \Longrightarrow (1)$ and $(3) \Longrightarrow (1)$ 
are easy consequences. \qed

\bigskip

Suppose that $\MM \prec  \NN$ and $I$ is a cut of $\MM$. We say that $\NN$ 
{\it fills} the cut $I$ if there is $b \in N$ such that 
$I = \{a \in M : \NN \models a < b\}$. In case $I = \omega$, we use the notation 
$\MM \prec ^{\omega} \NN$ to indicate that $\NN$ fills $\omega$. For us, the significance  of this definition 
is the following well known equivalence: If $\MM$ is a countable, nonstandard model, then $\ssy(\MM)$ is a jump ideal iff there is $\NN$ such that $\MM \prec^\omega \NN$ 
and $\ssy(\MM) = \ssy(\NN)$. 
See, for example,  \cite[Theorem~7.3.4]{ksbook}. 

\bigskip

The next lemma is  Lemma~3.11 of \cite{ks5}.

\bigskip

{\sc Lemma 2.4:} {\em Suppose that $\MM$ is a countable model and $I$ is a proper cut of $\MM$. If $I$ is not downward monotonically $\omega$-lofty, then there is a countable 
$\NN \succ_{\sf cf} \MM$ such that $\NN \not\succ^\omega \MM$ and $\sup^\NN(I)$ is an upward monotonically 
$\omega$-lofty cut of $\NN$.} 

\bigskip

{\it Proof}. We sketch the proof as it will be needed later on. The key notion that is used in the proof of Lemma~3.11 of \cite{ks5} 
is that of an  $\langle S,Q \rangle$-rich set (\cite[Def.\@ 3.6]{ks5}). We say that  $C$ is 
$\langle S,Q \rangle$-{\it rich} if $C,S,Q \in \Def(\MM)$, $S$ is bounded, and whenever 
$f : S \into Q$ is definable in $\MM$, then there is $c \in C$  such that $(c)_i = f(i)$ for every $i \in S$.

Suppose that $I$ is not downward monotonically $\omega$-lofty.

The key combinatorial fact that is used in the proof of Lemma~3.11 is the following:

\smallskip

{\sc Fact} 1: (\cite[Lemma 3.10]{ks5}) {\em Suppose that $r < \omega < s \in M$, $I < q \in M$, $C$ is 
$\langle [r,s), [0,q) \rangle$-rich, and $g : M \into M$ is definable. Then one of the following 
holds$:$
\begin{itemize}

\item[(1)] There are $n,t < \omega$ and $\{d_i : i \in [r,t)\} \subseteq I$ such that $r \leq t$ 
and $\{c \in C : g(c) \leq n$ and $(c)_i = d_i$ for $r \leq i < t\}$ is  $\langle [ t,s), [0,q)\rangle$-rich.

\item[(2)] There are $n,v,u \in M$ such that $n > \omega$, $\omega < v \leq s$, 
$I < u \leq q$ and $\{c \in C : g(c) > n\}$ is $\langle [r,v),[0,u) \rangle$-rich.

\end{itemize}}

\smallskip

The proof then proceeds as follows. Choose arbitrary  $s_0 = q_0 >I$, and let $r_0 = 0$ and $C_0 = M$. Obtain a decreasing sequence $C_0 \supseteq C_1 \supseteq C_2 \supseteq \cdots$ such that each $C_i$ is $\langle [r_i,s_i), [0,q_i) \rangle$-rich,  
where $r_i < \omega < s_i$ and $ I < q_i$. Also, there are $d_0,d_1, \ldots, d_{r_i-1} \in I$ such that whenever $c \in C_i$ and $j < r_i$, then $(c)_j = d_j$. Furthermore,  for each definable 
$g : M \into M$ there is an $i < \omega$ such that $C_{i+1}$ is obtained from $C_i$ in 
the obvious way using Fact~1. And finally, for each $d \in I$ there are $i,n < \omega$ 
such that $(c)_n = d$ for all $c \in C_i$. This sequence determines a type over 
$\MM$. Then, let $\NN$ be an elementary extension of $\MM$ generated by 
an element $c$ realizing this type. Clearly, $I = \{(c)_n : n <\omega\}$, 
so $\sup^\NN(I)$ is upward monotonically $\omega$-lofty. \qed

\bigskip

We get the following corollary (to be improved by Corollary 2.7 and again by Theorem~2.8).

\bigskip

{\sc Corollary 2.5:} {\em Suppose that $\MM$ is nonstandard and countable.  Then there is $\NN \succ_{\sf cf} \MM$ such that $\NN \not\succ^\omega \MM$ and 
$\NN$ has the $\omega$-property.}

\bigskip

{\it Proof}. Let $\MM_0 \succ_{\sf end} \MM$ be such that $\MM_0$ is countable and 
$M$ is not a downward 
monotonically $\omega$-lofty cut of $\MM_0$. (For example, let $\MM_0$ be a countable, conservative extension of $\MM$.) Apply Lemma~2.4 to get a countable $\NN_0 \succ_{\sf cf} \MM_0$ such that $\sup^{\NN_0}(M)$ is an upward monotonically $\omega$-lofty 
cut of $\NN_0$. Then the unique $\NN$ such that 
$\MM \preccurlyeq_{\sf cf} \NN \prec_{\sf end} \NN_0$ is as required. \qed

\bigskip

A construction, due to Paris,  of a countable $\NN$ that has the $\omega$-property but is 
not recursively saturated is presented in \cite[Theorem~3.2]{om1}. Another way of getting 
such a model $\NN$  is by use of the previous corollary. First, we need a definition taken from \cite[page~111]{ks4} 
(or see \cite[Notation 2.1(5)]{ks5}).  If $\MM$ is a model and $I \subseteq M$ is a cut, then we say that $I$ is 
{\it recursively definable} if there are an element $a \in M$ and two recursive 
sequences $t_0(x), t_1(x), t_2(x), \ldots$ and $t'_0(x), t_1'(x), t'_2(x), \ldots$ of Skolem terms such that 
$I = \inf\{t_i(a) : i < \omega\} = \sup\{t'_i(a) : i < \omega\}$. Obviously, no recursively saturated model has a recursively definable  cut. In fact, $\MM$ is recursively saturated iff it is tall and has no recursively definable cuts (\cite[Theorem~2.7(i)]{ks4}). There are prime models whose standard cuts are recursively definable 
(\cite[Theorem~2.3]{ks5}), and every completion of \pa\ has a finitely generated model 
whose standard cut is recursively definable 
(\cite[Coro.\@ 2.4]{ks5}). 

Now, start with a nonstandard countable model $\MM$ whose standard cut is recursively definable.  By Corollary~2.5, let $\NN \succ_{\sf cf} \MM$ be such that $\NN$ 
is countable, $\NN \not\succ^\omega \MM$ and $\NN$ 
 has the $\omega$-property. Since the standard cut of $\NN$ is  recursively definable, 
 $\NN$ is not 
recursively saturated.

\bigskip

We next show that  Corollary~2.5 can be somewhat improved.
To do so, we need the following variant of Lemma~2.4. 

\bigskip

{\sc Lemma 2.6:} {\em Suppose that $\MM$ is a countable model, $I$ is a proper cut 
of $\MM$ and $X \subseteq \omega$ is such that $X \not\in \ssy(\MM)$. 
If $I$ is not uniformly $\omega$-lofty, then there is a countable 
$\NN \succ_{\sf cf} \MM$ such that $\NN \not\succ^\omega \MM$, $\sup^\NN(I)$ is upward monotonically 
$\omega$-lofty and  $X \not\in \ssy(\NN)$.}

\bigskip

{\it Proof}. As in the proof of Lemma~2.4, we will obtain a decreasing sequence $C_0 \subseteq C_1 \supseteq C_2 \supseteq \cdots$ that  has all the properties  required of it in  that proof. But we also need an additional property that we get by interleaving into 
the construction some additional steps that are applications of the following fact.

\smallskip

{\sc Fact} 2: {\em  Suppose that $r < \omega < s \in M$, $I < q \in M$, $C$ is 
$\langle [r,s), [0,q) \rangle$-rich, and $g : M \into M$ is definable. Then there are 
$n,v,t,u \in M$ and $\{d_i: i \in [r,t)\} \subseteq I$ such that $n < \omega$, $r \leq t < \omega < v \leq s$, $I < u \leq q$ 
and $\{c \in C : (g(c))_n = 0 \Leftrightarrow n \not\in X $ and $(c)_i = d_i$ for $r \leq i < t\}$ is 
$\langle [t,v), [0,u) \rangle$-rich.}

\smallskip

We give a proof of this fact. Define two functions $f_0,f_1: M \times [r,s] \into [0,q]$ 
as follows. For each $n \in M$ and $t \in [r,s]$, let 
$f_0(n,t)$ be the largest $u \leq q$ such that 
$\{c \in C : (g(c))_n = 0\}$ is $\langle [r,t),[0,u) \rangle$-rich,  
and let $f_1(n,t)$ be the largest $u \leq q$ such that 
$\{c \in C : (g(c))_n > 0\}$ is $\langle [r,t),[0,u) \rangle$-rich. Both of these functions are well defined since both $\{c \in C : (g(c))_n = 0\}$ and $\{c \in C : (g(c))_n > 0\}$ are $\langle [r,t),[0,0) \rangle$-rich,  

Observe that for each $n \in M$ and $e \in \{0,1\}$, 
if $r \leq t_1 \leq t_2 \leq s$, then $f_e(n,t_1) \geq f_e(n,t_2)$. 

We consider three cases.

\smallskip

{\it Case 1}: There are $n < \omega$ and $v > \omega$ such that $f_0(n,v), f_1(n,v) > I$. 
We have two possibilities depending on whether or not $n \in X$. 

\smallskip

(1) $n \not\in X$: Let $u = f_0(n,v)$. Then, 
$\{c \in C : (g(c))_n = 0\}$ is $\langle [r,v), [0,u) \rangle$-rich.

 \smallskip
 
 (2)  $n \in X$: Let   $u =  f_1(n,v)$. Then, 
$\{c \in C : (g(c))_n >0\}$ is $\langle [r,v), [0,u) \rangle$-rich. 

\smallskip

{\it Case 2}: There are $n < \omega$ and $t < \omega$ such that $f_0(n,t), f_1(n,t) < I$. 
We again have two possibilities. 

\smallskip

(1) $n \not\in X$: Let $u = f_1(n,t) +1$. Since $\{c \in C : (g(c))_n > 0\}$ is not 
$\langle [r,t), [0,u) \rangle$-rich, there are $d_r,d_{r+1}, \ldots, d_{t-1} < u$ such that 
there is no $c \in C$ such that $(c)_i = d_i$ whenever $r \leq i < t$ and $(g(c))_n > 0$. 
Therefore, $\{c \in C : (g(c))_n = 0$ and $(c)_i = d_i$ for $r \leq i < t\}$ is 
$\langle [t, s), [0,u) \rangle$-rich. 

\smallskip

(2) $n \in X$: Let $u = f_0(n,t)+1$. Since $\{c \in C : (g(c))_n = 0\}$ is not 
$\langle [r,t), [0,u) \rangle$-rich, there are $d_r,d_{r+1}, \ldots, d_{t-1} < u$ such that 
there is no $c \in C$ such that $(c)_i = d_i$ whenever $r \leq i < t$ and $(g(c))_n = 0$. 
Therefore, $\{c \in C : (g(c))_n > 0$ and $(c)_i = d_i$ for $r \leq i < t\}$ is 
$\langle [t, s), [0,u) \rangle$-rich.

\smallskip

{\it Case 3}: Neither of the previous cases apply.  
Since $I$ is not uniformly $\omega$-lofty, there are $b,d$ such that $b < I < d \leq q$ 
and 
$[b,d) \cap \{f_0(n,t) : n< \omega, \ r \leq t < \omega\} = [b,d) \cap \{f_1(n,t) : n < \omega, \ r \leq t < \omega\} = \varnothing$.
By overspill, there is $u$ such that $\omega < u \leq s$ and 
$[b,d) \cap \{f_0(n,t) : n \leq u, \ r \leq t \leq u\} = [b,d) \cap \{f_1(n,t) : n \leq u, \ r \leq t \leq u\} = \varnothing$. Without loss of generality, we can take $u = s$. 

Define two functions $h_0,h_1 : [0,s) \into [0, q]$ as follows. For each $n < s$, let $h_0(n) = \max\{t \in [r,s] : f_0(n,t) \geq d\}$ and $h_1(n) = \max\{t \in [r,s]  : f_1(n,t) \geq d\}$. 
Both functions are well defined since  $f_0(n,r) = f_1(n,r) = q \geq d$. 

We claim, for all $n < \omega$, that $h_0(n) < \omega$ iff $h_1(n) > \omega$. 
For, if $h_0(n), h_1(n) > \omega$, then $v = \min(h_0(n),h_1(n))$ 
would put us into Case 1. And if $h_0(n),h_1(n) < \omega$, then 
$t = 1+ \max(h_0(n),h_1(n))$ would put us  into Case 2.

Let $Y = \{n < \omega : h_0(n) < h_1(n)\}$. Clearly, $Y \in \ssy(\MM)$, so $Y \neq X$.
Let $n < \omega$ be such that $n \in Y$ iff $n \not\in X$. We have two possibilities. 

\smallskip

(1) $n \not\in X$: Let $v = h_0(n) > \omega$ and  continue as in (2) of Case 1.  

\smallskip

(2) $n \in X$: Let $t = h_0(n)+1 < \omega$ and then continue as in (2) of Case 2. 

This completes the proof of Fact 2. 

We leave it to the reader to supply the missing details of the proof of the lemma. \qed

\bigskip

{\sc Corollary 2.7:} {\em Suppose that $\MM$ is nonstandard and countable, 
$X \subseteq \omega$ and $X \not\in \ssy(\MM)$.   Then there is $\NN \succ_{\sf cf} \MM$ such that $\NN \not\succ^\omega \MM$, $X \not\in \ssy(\NN)$ and 
$\NN$ has the $\omega$-property.}

\bigskip

{\it Proof}. Let $\MM_0$ be a countable conservative end extension of $\MM$. 
Then, one easily checks that $M$ is not a uniformly $\omega$-lofty cut of $\MM_0$. 
Now, proceed  as in the proof of Corollary~2.5 but invoking Lemma~2.6 instead of 
Lemma~2.4. \qed 

\bigskip

Both Corollaries~2.5 and 2.7 imply that  some countable, nonstandard  models $\MM$ that fail to be recursively saturated have cofinal extensions $\NN$ that have the $\omega$-property and also  fail 
to be recursively saturated. This can happen if  $\omega$ is a recursively definable cut 
(Corollary~2.5) or if $\Th(\MM) \not\in \ssy(\MM)$ (Corollary~2.7).  
We show next that, in fact, this is the case 
for all countable, nonstandard $\MM$ that are not recursively saturated.

\bigskip

{\sc Theorem 2.8}: {\em Suppose that $\MM$ is nonstandard, countable and not recursively saturated. 
   Then there is a countable $\NN \succ_{\sf cf} \MM$ such that 
   $\NN \not\succ^\omega \MM$,   
   $\NN$ has the $\omega$-property and $\NN$ is not recursively saturated.}

\bigskip

{\it Proof}. We can assume that $\MM$ is tall, as otherwise Corollary~2.5 does it. 
Since $\MM$ is tall and not recursively saturated, there are $a \in M$ and a recursive set $\Phi(x,y)$ of formulas such that $\Phi(x,y)$ includes $x< y$ and $\Phi(x,a)$  is 
finitely realizable in $\MM$ but  is 
omitted by $\MM$. The cofinal extension $\NN$ that we will obtain will fail to be recursively 
saturated by virtue of also omitting $\Phi(x,a)$. By \cite[Lemma~2.4]{ks4}, there is a recursively definable 
cut $I$  of $\MM$ such that for any cofinal extension $\NN$ of $\MM$, 
$\NN$ realizes $\Phi(x,a)$ iff $\NN$ fills $I$. So we will ignore 
trying to omit $\Phi(x,a)$ and focus instead on trying not to fill $I$.

By \cite[Lemma~3.17]{ks5}, there is a countable $\MM_0 \succcurlyeq_{\sf cf} \MM$ 
that fills neither $I$ nor $\omega$ and  has the property: 
whenever $\NN \succ \MM_0$, then $\NN$ fills $\omega$ iff $\NN$ fills $I^{\MM_0}$. 
Apply Corollary~2.5 to get $\NN \succ_{\sf cf} \MM_0$ such that $\NN \not\succ^{\omega} 
\MM_0$ and $\NN$ has the $\omega$-property. Then $\NN \not\succ^{\omega} 
\MM$ and $\NN$ does not fill $I^{\MM_0}$, so, as an extension of $\MM$, 
$\NN$ does not fill $I$. \qed

\bigskip

It is possible to get the natural common generalization of Lemma~2.6 and Theorem~2.8. 

Theorem~2.8 allows us to get  some characterizations of countable, recursively saturated models.

\bigskip

{\sc Corollary 2.9}: {\em Suppose that  $\MM$ is countable and tall. The following are equivalent$:$

\begin{itemize} 

\item[(1)]   $\MM$ is recursively saturated. 

\item[(2)] Every countable $\NN \succcurlyeq_{\sf cf} \MM$  is generated by a set of indiscernibles.

\item[(3)] Every countable $\NN \succcurlyeq_{\sf cf} \MM$ is isomorphic to 
some $\NN_0 \prec_{\sf end} \NN$.

\end{itemize}}

 \bigskip
 
 {\it Proof}. The implications   $(1) \Longrightarrow (2)$ and $(1) \Longrightarrow (3)$ are 
 consequences of Proposition~1.3.

We prove the converse implications.  Suppose that  $\MM$ is not recursively saturated. Let  $\NN$ be as in Theorem~2.8. By Lemma~2.2, $\MM$ is uniformly $\omega$-lofty. Then Lemma~2.3 implies that $\NN$ is not generated by a set of indiscernibles 
and there is no $\NN_0 \prec_{\sf end} \NN$ that is isomorphic to $\NN$. \qed

\bigskip

If the requirement in Corollary~2.9 that $\MM$ is tall is replaced with the weaker one that 
$\MM$ is nonstandard, then the equivalence $(1) \Longleftrightarrow (3)$ still holds. 
However, we don not know if $(1) \Longleftrightarrow (2)$ still holds.

We next give a characterization of those countable models that code their own theories.

\bigskip

{\sc Corollary 2.10}: {\em Suppose that $T$ is a completion of \pa\ and ${\mathfrak X}$ 
is a countable Scott set such that $\rep(T) \subseteq {\mathfrak X}$. 
The following are equivalent$:$

$(1)$ $T \not\in {\mathfrak X};$

$(2)$ Whenever $\NN \models T$ is  countable  and tall and 
$\ssy(\NN) = {\mathfrak X}$, then $\NN$ has a countable, cofinal extension $\NN_1 \succ_{\sf cf} \NN$  that  is not generated by a set of indiscernibles.}

\bigskip

{\it Proof}. $(1) \Longrightarrow (2)$: Suppose that $T \not\in {\mathfrak X}$. Consider 
an arbitrary countable  and tall $\NN \models T$  such that $\ssy(\NN) = {\mathfrak X}$. 
By Corollary~2.7, let 
$\NN_1 \succ_{\sf cf} \NN$ be a countable model that has the $\omega$-property 
such that $T \not\in \ssy(\NN_1)$. By Lemma~2.2, $\NN_1$ is uniformly $\omega$-lofty. 
Since $T = \Th(\NN_1) \not\in \ssy(\NN_1)$, then $\NN_1$ is not recursively saturated, 
so $\NN_1$ is not generated by a set of indiscernibles according to  Lemma~2.3. 

\smallskip

$(2) \Longrightarrow (1)$: Suppose that $T \in {\mathfrak X}$. By Theorem~1.1, we can let 
$\NN$ be a countable, recursively saturated model of $T$ such that $\ssy(\NN) = {\mathfrak X}$. By Propositions 1.3(2) and (4),  every countable $\NN_1 \succ_{\sf cf} \NN$ is generated by 
a set of indiscernibles. \qed

\bigskip


{\bf \S3. Substructure Lattices.}  It was proved in \cite[Theorem~5.1]{ks2} that if 
$\MM, \NN$ are arithmetically saturated models and $\Lt(\MM) \cong \Lt(\NN)$,
then $\ssy(\MM) \cong \ssy(\NN)$. In this section we will take a closer look at the proof of that result 
in order to obtain a refinement of it stated below as Theorem~3.2. This theorem follows 
immediately from Lemma~3.1, which is the main result of this section and  will be 
applied later on in proving the main results of this paper.

If $L = (L, \vee, \wedge)$ is a lattice, we will sometimes think of it as a partially ordered 
set $(L, \leq)$, where $x \leq y$ iff $x = x \wedge y$ iff $y = x \vee y$. We let $0_L$ be the  least element of $L$ if there is one, and we  let $1_L$ be the largest element 
if there is one. If $n < \omega$, then ${\mathbf n}$ is the lattice that is a chain having exactly $n$ elements. We let ${\bf B}_2$ be the 4-element Boolean lattice; that is, 
${\bf B}_2$ is the unique 4-element lattice that is not ${\mathbf 4}$. If $K,L$ are two lattices and both $1_K, 0_L$ exist, 
then their {\it linear sum} $K \oplus L$ is the lattice that is the disjoint union of $K$ and $L$ 
(except that we set $1_K = 0_L$) such that both $K,L$ are sublattices of $K \oplus L$ and 
$x \leq y$ whenever $x \in K$ and $y \in L$. (For example, ${\mathbf 2} \oplus {\mathbf 3} = {\mathbf 4}$.)

We next define some more lattices that are the same ones defined in \cite[\S5]{ks2}.
If $X \subseteq n < \omega$, we define the lattice ${\mathcal D}_0(X,n)$ by recursion on $n$ as follows.
Let 
$${\mathcal D}_0(\varnothing, 0) = {\mathbf 1},$$ 
 and then if $X \subseteq n+1$, 
let
\begin{displaymath}
{\mathcal D}_0(X,n+1) = \left\{ 
     \begin{array}{ll}
     {\mathcal D}_0(X \cap n,n) \oplus {\mathbf 2} & {\mbox{ if }} \ n  \in X \\ 
       {\mathcal D}_0(X \cap n,n) \oplus {\bf B}_2 & {\mbox{ if }} \  n \not\in X.
       \end{array} \right.
       \end{displaymath}
 Note that $|{\mathcal D}_0(X,n)| = 3n+1-2|X|$. 
 (Roughly, ${\mathcal D}_0(X,n)$ is a stack of $n$  lattices, the $i$-th one in the stack 
 is ${\mathbf 2}$ iff $i \in X$ and is ${\bf B}_2$ iff $i \not\in X$.)
  If $X \subseteq n+1$, then 
 ${\mathcal D}_0(X \cap n, n)$ is an ideal of ${\mathcal D}_0(X,n+1)$. Next, if $X \subseteq n$, let 
 ${\mathcal D}(X,n) = {\mathcal D}_0(X,n) \oplus {\mathbf 2}$. If $X \subseteq n+1$, we consider 
 that ${\mathcal D}(X \cap n, n) \subseteq {\mathcal D}(X,n+1)$ by setting  $1_{{\mathcal D}(X \cap n,n)} = 1_{{\mathcal D}(X,n+1)}$. 
  
If $X \subseteq \omega$, we define 
${\mathcal D}(X)  = \bigcup_{n<\omega}({\mathcal D}(X \cap n,n)$,
and then let ${\mathcal D}'(X)$ be the lattice obtained from ${\mathcal D}(X)$ by adding one more element 
that is less than $1_{{\mathcal D}(X)}$ but greater than every other element of ${\mathcal D}(X)$. Thus, each ${\mathcal D}(X \cap n,n)$ is a sublattice of 
${\mathcal D}(X)$ which, in turn, is a sublattice of ${\mathcal D}'(X)$. 
In fact, ${\mathcal D}'(X)$ is the completion of ${\mathcal D}(X)$, and ${\mathcal D}(X)$ is the sublattice of ${\mathcal D}'(X)$ 
consisting of the compact elements in ${\mathcal D}'(X)$.

\bigskip

{\sc Lemma 3.1:} {\em Suppose that $\MM$ is  recursively saturated and  $X \subseteq \omega$. The following are equivalent.

\begin{itemize}

\item [(1)] There is $\MM_1 \prec \MM$ such that $\Lt(\MM_1) \cong {\mathcal D}'(X)$.

\item [(2)] There is $Y \in \ssy(\MM)$ such that $X \leq_T Y'$.

\end{itemize}}

\bigskip

{\it Proof.} $(1) \Longrightarrow (2)$: Let $\MM_1 \prec \MM$ be such that $\Lt(\MM_1) \cong {\mathcal D}'(X)$. Since $1_{{\mathcal D}'(X)}$ exists and has a unique immediate predecessor, it must be that $\MM_1$ is finitely generated. Let $\MM_1 = \scl(a)$, and $Y = \tp(a)$. Since, 
$\MM$ is recursively saturated, $Y \in \ssy(\MM)$. 

We easily see that 
$X$ is $\Delta^0_2$ in $Y$. First,  $n \in X$ iff  there is $Z \subseteq n+1$ such that $n \in Z$ and there are Skolem terms $t_r(x)$ for each $r \in {\mathcal D}_0(Z,n)$ (where we let $m = 1_{{\mathcal D}_0(z,n)}$ and also
Skolem terms $t_{rs}(y)$ for each $r \leq s \in L(Z,n)$ such that each formula 
$t_{rs}(t_s(x)) = t_r(x)$ is in $\tp(a)$ and whenever $t(y)$  is a Skolem term, 
then for each $s \in {\mathcal D}_0(Z,n)$ there is  a unique $r \leq s$ such that 
$t(t_s(x)) = t_r(x)$ is in $\tp(a)$. Similarly, $n \not\in X$ iff there is $Z \subseteq n +1$ such that 
$n \not\in Z$ and all the same  conditions hold for this $Z$.   

\bigskip

$(2) \Longrightarrow (1)$: This part of the proof relies on the technology for constructing models 
with a prescribed substructure lattice as presented in \cite[Chap.\@ 4.5]{ksbook}. 
Since the application of this technology is quite routine, we will present just a sketch of the proof.

Let $\MM_0$ be the prime elementary submodel of $\MM$. For some finite lattices $L$, 
we define by recursion a representation $\alpha_L : L \into \eq(A)$, which we will call a
{\it regular} representation. Each regular representation that we define will be 
 definable in $\MM_0$, and we will define it only up to $\MM_0$-definable isomorphism.

\begin{itemize}

\item If $L = {\mathbf 2}$, then $\alpha_L$ is {\it regular} if $A = M_0$.

\item If $L = {\bf B}_2$, then assume that $b_1,b_2$ are the two atoms of ${\bf B}_2$. 
Then, $\alpha_L$ is {\it regular} if $A = \{\langle x_1,x_2 \rangle \in M_0^2 : x_1 < x_2\}$, 
and whenever $\langle x_1,x_2 \rangle, \langle y_1,y_2 \rangle \in A$ and $e \in \{1,2\}$, then 
$\big\langle \langle x_1,x_2 \rangle,$ $ \langle y_1,y_2 \rangle\big\rangle \in \alpha(b_e)$ iff 
$x_e = y_e$. 

\item If  $\alpha_0 : L_0 \into \eq(A_0)$ and $\alpha_1:L_1 \into \eq(A_1)$ are regular 
and $L = L_0 \oplus L_1$,  then $\alpha : L \into \eq(A)$ is {\it regular} if
$A = A_0 \times A_1$ and whenever 
$\langle x_0,x_1 \rangle, \langle y_0,y_1 \rangle \in A$ and 
$r \in L_0 \oplus L_1$, then 
$\big\langle \langle x_0,x_1 \rangle, \langle y_0,y_1 \rangle\big\rangle \in \alpha(r)$ 
iff one of the following holds:

\smallskip

(0) $r \in L_0$ and $\langle x_0,y_0 \rangle \in \alpha_0(r)$,

 (1) $r \in L_1$, $x_0 = y_0$ and  $\langle x_1,y_1 \rangle \in \alpha_1(r)$.

\smallskip

\end{itemize}

\noindent Whenever $L = {\mathcal D}_0(X,n)$ or $L = {\mathcal D}(X,n)$, the regular representation is well defined and  unique (up to $\MM$-definable isomorphism). 

We observe some facts about regular representations of the ${\mathcal D}(X,n)$'s.  Suppose that 
$\alpha : {\mathcal D}(X,n) \into \eq(A)$ is a regular representation. 

\begin{itemize}

\item[(1)] If $m < n$, then $\alpha \harp {\mathcal D}(X \cap m,m)$ is regular.

\item[(2)] If $n < m < \omega$, $Y \subseteq m$ and $X = Y \cap n$, then there a regular 
representation $\beta : {\mathcal D}(Y,m) \into \eq(A)$ such that $\alpha = \beta \harp {\mathcal D}(X,n)$.

\item[(3)]  If $\Theta \in \eq(A)$ is definable in $\MM_0$, 
then there is $B \subseteq A$ such that $\alpha|B$ is regular and there is $r \in {\mathcal D}(X,n)$ 
such that $\alpha(r) \cap B^2 = \Theta \cap B^2$. 

\end{itemize}
The proofs of (1) -- (3) will be omitted. The reader can consult with \cite[\S5]{ks2} 
where some analogous statements are proved.

Each of (1) -- (3) is effective. This is trivial for (1). For (2), this means that given $n < m < \omega$, 
$Y \subseteq m$ and a formula that defines $\alpha$, then a formula that 
defines $\beta$ can be effectively obtained. For (3), this means that given $n < \omega$ 
and formulas that define $\alpha$ and $\Theta$, then a formula defining $B$ can be effectively obtained. 

Since $\Th(\MM) \in \ssy(\MM)$, we can assume without loss of generality that $\Th(\MM) \leq_T Y$.

Let $\langle X_n : n < \omega \rangle$ be a sequence of finite sets that is recursive in $Y'$ 
such that $\lim_n X_n = X$ in the sense that for all $i < \omega$,
 there is $m < \omega$ such that whenever $m \leq n < \omega$, 
then $i \in X$ iff $i \in X_n$. We can also require that there is a sequence 
$0 = k_0,k_1,k_2, \ldots$ such that whenever $n < \omega$, then 
$X_n \subseteq k_n < \omega$ and either $k_{n+1} = k_n+1$ and 
$ X_n = X_{n+1} \cap n$ or else $k_{n+1} \leq k_n$ and $X_{n+1} = X_n \cap n+1$. 

Let $\theta_0(x,y), \theta_1(x,y), \theta_2(x,y), \ldots$ be a recursive list of 
$2$-ary formulas in the language of \pa\ so that each definable equivalence relation 
$\Theta \subseteq M_0^2$ is defined by infinitely many of the formulas. We construct a sequence of regular representations $\langle \alpha_n : {\mathcal D}(X_n,k_n) \into A_n \rangle$. (More precisely, we construct a sequence  
$\varphi_0,\varphi_1, \varphi_2, \ldots $ of formulas such that $\varphi_n$ defines 
$\alpha_n$ in $\MM_0$, and this sequence should be recursive in $Y'$.) 

Let $\alpha_0 : {\mathbf 2} \into \eq(M_0)$ be a regular representation. 
Now suppose that we $\alpha_n$. We will effectively obtain $\alpha_{n+1}$.
If $k_{n+1} = k_n+1$, apply (2) to get a regular $\alpha : {\mathcal D}(X_{n+1},k_{n+1}$, 
and if $k_{n+1} \leq k_n$, then let 
$\alpha = \alpha_n \harp {\mathcal D}(X_{n+1},k_{n+1})$. Since $\alpha$ is regular,
we apply (3) with $\Theta$ being the equivalence relation defined bt $\theta_n$ to get 
$\alpha_{n+1} = \alpha|B$.

Thus, for each $n < \omega$, we have  $\alpha_n : L(X_n,k_n) \into \eq(A_n)$. 
The sequence $A_0 \supseteq A_1 \supseteq A_2 \supseteq$ is recursive and determines 
a complete type $p(x)$. Let $a$ realize this type in $\MM$, and then let 
$\MM_1 = \scl(a)$. 

One then checks that $\Lt(\MM_1) \cong {\mathcal D}'(X)$. \qed

\bigskip

The previous lemma easily implies the following theorem.

\bigskip

{\sc Theorem 3.2:} {\em If $\MM, \NN$ are recursively saturated models 
and $\Lt(\MM)$ $ \cong \Lt(\NN)$, then for each $X \in \ssy(\MM)$ there is $Y \in \ssy(\NN)$ 
such that $X' \equiv_T Y'$.} \qed

\bigskip

If $\MM$ is recursively saturated, then ${\bf M}_3$ is an ideal of $\Lt(\MM)$ iff $\MM$ 
is not a model of $\ta$. Can Theorem~3.2 be improved if neither (or both) are models of \ta?

\bigskip


{\bf \S4.\@ The Proof of Theorem 4.}  Lascar's Theorem shows (in the terminology of \cite{ns}) that open subgroups are recognizable  for the class of countable, arithmetically saturated models. This means: 
if $\MM, \NN$ are countable, arithmetically saturated models, $\alpha : \aut(\MM) \into \aut(\NN)$ 
is an isomorphism and $H \leq \aut(\MM)$, then $H$ is an open subgroup of $\aut(\MM)$ iff 
$\alpha[H]$ is an open subgroup of $\aut(\NN)$. 
This consequence of Lascar's Theorem was improved in   \cite[Coro.\@ 3.14]{ns} where it was shown that basic open 
subgroups are recognizable.{\footnote{Henceforth, the term ``recognizable'' will be used in 
a rather informal way, and it is be understood as meaning ``recognizable for the class of 
arithmetically saturated models''.}} 
In particular, 
if $\MM, \NN$ are countable, arithmetically saturated models, $\alpha : \aut(\MM) \into \aut(\NN)$ 
is an isomorphism and $H \leq \aut(\MM)$, then $H$ is the stabilizer of a finite set iff 
$\alpha[H]$ (as a subgroup of $\aut(\NN)$) is the stabilizer of a finite set. This allows us to define 
the function $\wt{\alpha}: \Lt_0(\MM) \into \Lt_0(\NN)$ as follows: If $\MM_1 \in \Lt_0(\MM)$, then 
  $\wt{\alpha}(\MM_1)$ is that unique $\NN_1 \in \Lt_0(\NN)$ such that 
  $ \aut(\NN)_{(N_1)} = \alpha[\aut(\MM)_{(M_1)}] $. We easily see that $\wt{\alpha}$ is an isomorphism
  from the semilattice $\Lt_0(\MM)$ onto the semilattice $\Lt_0(\NN)$. Thus, it extends uniquely 
  to an isomorphism of the lattices $\Lt(\MM)$ and $\Lt(\NN)$. We denote this extension also by 
  $\wt{\alpha}$. Thus, we have the following lemma.

  \bigskip

{\sc Lemma 4.1}: {\em  Suppose that $\MM, \NN$ are countable, arithmetically saturated models 
and that $\alpha : \aut(\MM) \into \aut(\NN)$ is an isomorphism. 
Then $\wt{\alpha} : \Lt(\MM) \into \Lt(\NN)$ is an isomorphism.}

\bigskip

We will typically invoke this lemma without referencing it.

 The isomorphism $\wt{\alpha}$ is implicit in \cite[Coro.\@ 3.15]{ns}. 
 Observe that if $a \in M$ and $g \in \aut(\MM)$, then $g \in \aut(\MM)_a$ iff $\alpha(g) \in \wt{\alpha}(\aut(\MM)_a)$. 
 
The map $\alpha \mapsto \wt{\alpha}$ is functorial in the sense that if $\MM_1, \MM_2, \MM_3$ 
are countable, arithmetically saturated models and  $\alpha : \aut(\MM_1) \into \aut(\MM_2)$, $\beta : \aut(\MM_2) \into \aut(\MM_3)$
are isomorphisms, then $\wt{\beta\alpha} = \wt{\beta}\wt{\alpha}$ and 
$\wt{\alpha^{-1}} = \wt{\alpha}^{-1}$.

\bigskip

We can now  comment about the general strategy that is used in this section. Of course, the ultimate goal is to prove Theorem~4. It will be seen that Theorem~4 follows 
almost immediately from Lemma~3.1 and the fact, to be proved as Lemma~4.15, that 
recursively saturated structures are recognizable. To get that conclusion, we will prove 
that recursive saturation can be characterized in terms of properties already shown to be recognizable. But to get the recognizability of these other properties, we will show that they too are characterizable in terms of other properties that were previously shown to be recognizable. And so on. Thus, we will build a catalogue of recognizable properties and show  that recursive saturation is in this catalogue. 

We begin by showing that isomorphism is recognizable

\bigskip

{\sc Lemma 4.2}: {\em Suppose that $\MM, \NN$ are countable, arithmetically saturated models 
and that $\alpha : \aut(\MM) \into \aut(\NN)$ is an isomorphism. 
If $\MM_1, \MM_2 \preccurlyeq \MM$ and $\MM_1 \cong \MM_2$, then $\wt{\alpha}(\MM_1) \cong \wt{\alpha}(\MM_2)$.} 

\bigskip

{\it Proof}. Let $G = \aut(\MM)$. Suppose that $\MM_1, \MM_2 \preccurlyeq \MM$ and that $\MM_1 \cong \MM_2$. Let 
$\NN_1 = \wt{\alpha}(\MM_1)$ and $\NN_2 = \wt{\alpha}(\MM_2)$. We wish to show that 
 $\NN_1 \cong \NN_2$.

First, suppose that  $\MM_1$ is finitely generated.   Then, $\MM_2$ is finitely generated, and 
$G_{(M_1)}$ and $G_{(M_2)}$ are conjugate 
subgroups, so let $g \in G$ be such that 
$gG_{(M_1)}g^{-1} = G_{(M_2)}$.   Applying $\alpha$ yields 
$$\alpha(g)\aut(\NN)_{(N_1)}(\alpha(g))^{-1} = \aut(\NN)_{(N_2)},$$ so that 
$\aut(\NN)_{(N_1)}$, $\aut(\NN)_{(N_2)}$ are conjugate 
subgroups of $\aut(\NN)$. Therefore,   $\NN_1 \cong \NN_2$. 

We have just proved the lemma in the case $\MM_1, \MM_2$ are finitely generated. 
Observe that if $f \in G$, then $f \harp M_1$ is an isomorphism from $\MM_1$ onto 
$\MM_2$ iff $f G_{(M_1)} f^{-1} = G_{(M_2)}$.

Next, suppose that $\MM_1 \preccurlyeq \MM$ is not finitely generated, and that 
$\varphi : \MM_1 \into \MM_2$ is an isomorphism.  Let $\MM_{1,0} \prec \MM_{1,1} \prec \MM_{1,2} \prec \cdots$ be a sequence of 
finitely generated, elementary substructures of $\MM_1$ such that 
$\MM_1 = \bigcup_{i< \omega} \MM_{1,i}$.  For each $i < \omega$, let $\varphi_i = \varphi \harp M_{0,i}$ and 
$\MM_{2,i} = \varphi_i[\MM_{1,i}]$. Then, $\MM_2 = \bigcup_{i< \omega} \MM_{2,i}$. 
Let $f_i \in G$ be such that $f_i \supseteq \varphi_i$. 
Since $\varphi_0 \subseteq \varphi_1 \subseteq \varphi_2 \subseteq \cdots$,
we have that whenever $i \leq j < \omega$, then  
$f_j G_{(M_{0,i})}f^{-1} = G_{(M_{1,i})}$.

 Let $\NN_{1,i} = \wt{\alpha}(\MM_{1,i})$ and $\NN_{2,i} = \wt{\alpha}(\MM_{2,i}$. 
 Then,  $\NN_1 = \bigcup_{i< \omega} \NN_{1,i}$ and $\NN_2 = \bigcup_{i< \omega} \NN_{2,i}$.  By the  first part of this proof, each $\NN_{0,i} \cong \NN_{1,i}$. Let $\theta_i : \NN_{0,i} \into \NN_{1,i}$ be the unique isomorphism. Let $g_i = \alpha(f_i)$.  
 Then, whenever $i \leq j < \omega$, then  
$g_j \aut(\NN)_{(N_{0,i})}g_j^{-1} = \aut(\NN)_{(N_{1,i})}$. Thus,
if  $i \leq j < \omega$, then $g_j\harp N_{0,i} : \NN_{0,i} \into \NN_{1,i}$ is an isomorphism.
Hence, 
$\bigcup_{i<\omega}g_i \harp N_{0,i}$ is an isomorphism from $\NN_0$ onto $\NN_1$. 
  \qed

\bigskip

The next lemma says that both cofinal extensions and end extensions are recognizable 
and also that both tall models and short models are recognizable.

\bigskip 

{\sc Lemma 4.3}: {\em Suppose that $\MM, \NN$ are countable, arithmetically saturated models, $\alpha : \aut(\MM) \into \aut(\NN)$ is an isomorphism and 
 $\MM_1, \MM_2$ $ \preccurlyeq \MM$.

\begin{itemize}

\item[(1)]   If $\MM_1 \prec_{\sf cf} \MM_2$, then $\wt{\alpha}(\MM_1) \prec_{\sf cf} \wt{\alpha}(\MM_2)$.

\item[(2)]   If $\MM_1 \prec_{\sf end} \MM_2$, then $\wt{\alpha}(\MM_1) \prec_{\sf end} \wt{\alpha}(\MM_2)$.

\item[(3)] If $\MM_1$ is tall, then $\wt{\alpha}(\MM_1)$ is tall.

\item[(4)] If $\MM_1$ is short, then $\wt{\alpha}(\MM_1)$ is short.

\end{itemize}}

\bigskip

{\it Proof}. First, note that the special case of (2) in which $\MM_2 = \MM$ is  one of the parts of  \cite[Prop.\@ 2.1]{ns} and is also \cite[Coro.\@  9.4.10]{ksbook}, where a proof is given. For us, an easier proof comes from an application of  Kaye's Theorem, 
or, more precisely, the $*$-version of Kaye's Theorem that applies 
to structures of the form $(\MM,a)$, where $a \in M$.  Observe that  $H < \aut(\MM)$ is the pointwise stabilizer of a short  elementary cut iff for some $a \in M$, 
 $H$ is the smallest, nontrivial, closed normal subgroup of $\aut(\MM)_a = \aut(M,a)$. Then note that every elementary cut is the union of a set of short elementary cuts, 
and, conversely, the union of any set of short elementary cuts is an elementary cut.

We now prove (1) -- (4). 

\smallskip

(1):\@ 
 This is due to the 
 following equivalence: 
$\MM_1 \prec_{\sf cf} \MM_2$ iff $\MM_1 \prec \MM_2$ and whenever 
$\MM_1 \preccurlyeq \MM_3 \prec_{\sf end} \MM$, then $\MM_2 \preccurlyeq \MM_3$.

\smallskip

(2):\@ This is a consequence of (1), and the following equivalence: 
$\MM_1 \prec_{\sf end} \MM_2$ iff $\MM_1 \prec \MM_2$ and there is no $\MM_3$ 
such that $\MM_1 \prec_{\sf cf}  \MM_3 \preccurlyeq \MM_2$. 

\smallskip

(3):\@ This is a consequence of (2) since $\MM_1$ is tall iff there is a 
sequence $\MM_{1,0} \prec_{\sf end} \MM_{1,1} \prec_{\sf end}  \MM_{1,2}
\prec_{\sf end} \cdots$ such that $\MM_1 = \bigcup_{i< \omega} \MM_{1,i}$.

\smallskip

(3):\@ $\MM_1$ is short iff $\MM_1$ is not tall.
\qed

\bigskip 

Suppose that $\MM$ is arithmetically saturated. The interstices and interstitial gaps 
of $\MM$ were first defined and studied in \cite{bk}. The {\it least interstice} of $\MM$, 
denoted by $\Omega_\omega$, is the set  $\{x \in M : \omega < x < a$ for all $a \in \scl(\varnothing) \backslash \omega\}$.  The arithmetic saturation 
of $\MM$ entails that $\Omega_\omega \neq \varnothing$. Observe that $\Omega_\omega \cup \hspace{2pt}   \omega$ is the smallest nonstandard invariant cut of $\MM$. We partition $\Omega_\omega$ into  convex sets, called interstitial gaps (or {\it igaps} for short) 
as follows.
 First, let ${\mathcal F}$ be the set of $\varnothing$-definable functions
$f : M \into M$ such that whenever $x \leq y < \omega$, then $x \leq f(x) \leq f(y) < \omega$. Then, for each $a \in \Omega_\omega$, define the
 {\it igap around} $a$ to be the set
$$\igap(a) = \{b \in \Omega_\omega : a \leq f(b) {\mbox{ and }} b \leq f(a) {\mbox{ for some }} f \in F^{\MM}\}.$$
The set of igaps is linearly ordered with order type of the rationals. It is a routine exercise 
in recursive saturation to show that whenever $\gamma_1 < \gamma_2$ are igaps and 
$a \in \Omega_\omega$, then there is $g \in \aut(\MM)$ such that $\gamma_1 < g(a) < \gamma_2$. 
A cut $I \subseteq M$ is an  {\it icut} if $I \subseteq  \omega \cup \Omega_\omega$ and 
whenever $\gamma$ is an  igap such that $\gamma \cap I \neq \varnothing$, 
then $\gamma \subseteq  I $. 

Kaye's Theorem implies that $\aut(\MM)_{(\Omega_\omega)}$ is the largest, closed 
proper normal subgroup of $\aut(\MM)$.

The following lemma  improves Proposition 4.2 of \cite{ns} by  eliminating the hypothesis that both $\MM$ and $\NN$ are 2-Ramsey.{\footnote{According to \cite{ns}, a model $\MM$ 
is $n$-Ramsey iff $({\mathbb N}, \rep(\Th(\MM))) \models {\sf RT}^n_2$.}}

\bigskip

{\sc Lemma 4.4}: {\em Suppose that $\MM, \NN$ are countable, arithmetically saturated models 
and that $\alpha : \aut(\MM) \into \aut(\NN)$ is an isomorphism. 
Suppose that $H_1 \leq \aut(\MM)$ and $H_2 = \alpha[H_1] \leq \aut(\NN)$. 
\begin{itemize} 
\item[(a)] If $H_1$ is a pointwise stabilizer of an icut, then so is $H_2$.
\item[(b)] If $H_1$ is a setwise stabilizer of an icut, then so is  $H_2$.
\item[(c)] If $H_1$ is a pointwise stabilizer of an igap, then   so is $H_2$.
\item[(d)] If  $H_1$ is a setwise stabilizer of an igap, then  so is $H_2$.
\end{itemize}}

\bigskip

{\it Proof}. We will prove (a). Then parts (b) -- (d) will follow just as in the proof of \cite[Prop.\@ 4.2]{ns}. 

Let $G = \aut(\MM)$. The concepts of least interstice, igap and icut extend naturally to models 
$(\MM,a)$, where $a \in M$.  

Suppose that $a \in M$. Let $J(a) = \sup((\omega \cup \Omega_\omega) \cap \scl(a))$ 
and then let $N(a) = G_{(J(a))} \cap G_a$.  It is easily checked that $J(a)$ is an icut. 
By the $*$-version of Kaye's Theorem (applied to $(\MM,a)$), 
$N(a)$ is  
  the smallest closed  subgroup $H$ such that $G_{(\Omega_\omega)} \cap G_a < H \trianglelefteq G_a$. 
  
 We will say that two groups $N(a)$ and $N(b)$ are {\it equivalent} if there is $c \in M$ 
 such that $N(a) \cap G_c = N(b) \cap G_c$. 
 
 \smallskip
 
 {\sf Claim}: {\em If $a,b \in M$, then  $N(a)$ and $N(b)$ are equivalent iff $J(a) = J(b)$.}
 
 \smallskip
 
 To prove the claim, consider $a,b \in M$. 
  If $J(a) = J(b)$, then $c = \langle a,b \rangle$ is such that 
  $N(a) \cap G_c = N(b) \cap G_c$. 
  
 For the converse, suppose that $J(b) \backslash J(a) \neq \varnothing$. Then, by recursive saturation, for any $c \in  M$,
there is $f \in G_{(J(a)} \cap G_c$ that moves some $d \in J(b)$. This proves the claim.
 
 \smallskip 
 
 It  follows from the arithmetic saturation of $\MM$ that for any $c \in M  \backslash \scl(J(a))$, there is an $f \in G_{(J(a))}$ 
 such that $f(c) \neq c$. Thus, we have that 
 $G_{(J(a))}$ is the closure of  $\bigcup\{N(b) : N(b)$ is equivalent to $N(a)\}$. Thus, for any $a \in M$ there is $b \in N$ such that $\wt{\alpha}(\scl(N(a)) = \scl(N(b))$.

 Now suppose that $I$ is an arbitrary icut. Either there is an igap such that 
 $I = \sup(\gamma)$ or there is not. If there is no such $\gamma$, then 
 $I$ is the union of all those $J(a)$ such that $J(a) \subseteq I$, and if there is such 
 a $\gamma$, then $I$ is the intersection of all those $J(a)$ such that $J(a) \supseteq I$. 
 Thus, $H$ is the pointwise stabilizer of an icut iff it is either the union or the intersection 
 of a set of subgroups of the form $G_{(J(a))}$, from which  
  (a) easily follows. 
 \qed

 \bigskip

In the next lemma, we formally record that certain  kinds of models  are recognizable. 
In particular, those described in 
(2) and (3) of Lemma~2.3 are. The recognizability of those models 
in (1) of Lemma~2.3 will be shown in Lemma~4.15.

\bigskip

  {\sc Lemma 4.5}:   {\em Suppose that $\MM, \NN$ are countable, arithmetically saturated models 
and that $\alpha : \aut(\MM) \into \aut(\NN)$ is an isomorphism. 
Suppose that $\MM_0 \prec \MM$.

\begin{itemize}

\item[(1)] If $\MM_0$ is homogeneous, then $\wt{\alpha}(\MM_0)$ is homogeneous.

\item[(2)] If $\MM_0$ is generated by a set of indiscernibles, then so is 
  $\wt{\alpha}(\MM_0)$.

\item[(3)] If $\MM_0$ is isomorphic to some $\MM_1 \prec_{\sf end} \MM_0$, then 
 $\wt {\alpha}(\MM_0)$ is isomorphic to some $\NN_1 \prec_{\sf end} \wt{\alpha}(\MM_0)$.

\end{itemize}}

\bigskip

{\it Proof}. The proofs are straightforward. We only note that (1) follows from the 
characterization: 
 $\MM_0$ is homogeneous iff whenever 
$a_0,a_1,b_0 \in M_0$ and $\tp(a_0) = \tp(b_0)$, then there is $b_1 \in M_0$ such that 
$\tp(a_0,a_1) = \tp(b_0,b_1)$. \qed

  \bigskip
  
  One sort of problem that arises when working inside an arithmetically saturated 
  model $\MM$ is that some familiar constructions cannot be carried out. 
  For an example, take the MacDowell-Specker Theorem: 
  For every $\MM_0$ there is $\MM_1$ such that $\MM_0 \prec_{\sf end } \MM_1$. 
  However, even if $\MM$ is arithmetically saturated, there is $\MM_0 \prec \MM$ for which there  is no $\MM_1$  
  such that $\MM_0 \prec_{\sf end } \MM_1 \preccurlyeq \MM$. One way around this problem is by restricting attention to just the coded elementary substructures of $\MM$. 
  We next make the appropriate definitions.
  
  If $\MM$ is a model, then a subset $A \subseteq M$ is  {\it coded} ({\it in} $\MM$) if there is 
$a \in M$ such that $A  = \{(a)_i : i < \omega\}$. In particular, $\ssy(\MM)$ is the set of coded subsets of $\omega$. A model $\MM_0$ is {\it coded} in $\MM$ if $\MM_0 \prec \MM$ and $M_0$ is coded. If $\MM$ is recursively saturated and $A \subseteq M$ is coded, then $\scl(A)$ is coded. Clearly, every finitely generated 
elementary substructure of a recursively saturated model is coded. 
If $\MM$ is recursively saturated and $\MM_0 \prec \MM$ is coded, then there is a coded 
$\MM_1 \prec \MM$ such that $\MM_0 \prec_{\sf end} \MM_1$. The following 
proposition shows that even more is true.

\bigskip

{\sc Proposition 4.6}: {\em Suppose that  $\MM$ is recursively saturated and $\MM_0 \prec \MM$ is coded.  Then there is a coded, tall  
$\MM_1 \prec \MM$ such that $\MM_0 \prec_{\sf end} \MM_1$ and $\MM_1$ is a conservative extension of $\MM_0$. }

\bigskip

{\it Proof}. Suppose that $\MM_0$ is coded. Let $a \in M$  realize a minimal type 
in $\MM$. By recursive saturation, there is $b \in M$ such that for each $i < \omega$, 
$M_0 < (b)_i < (b)_{i+1}$ and $\tp((b)_i) = \tp(a)$. 
Then $\MM_1 = \scl(M_0 \cup \{(b)_i : i < \omega\})$ is coded and tall and is a conservative 
extension of $\MM_0$. 
 \qed

\bigskip

Recall the notion of a superminimal extension as in \cite[\S2.1.2]{ksbook}. 
If $\MM_1 \prec \MM_2$, then $\MM_2$ is a {\it superminimal} extension 
of $\MM_1$ if $\MM_3 \preccurlyeq \MM_1$ for every $\MM_3 \prec \MM_2$. 
It is clear that superminimal extensions are recognizable.

\bigskip

{\sc Lemma 4.7}: {\em Suppose that $\MM$ is arithmetically saturated and $\MM_1 \preccurlyeq \MM$. Then $\MM_1$ is coded iff $\MM_1$ has a  superminimal 
end extension $\MM_2 \prec \MM$.}

 \bigskip
 
 {\it Proof}. $(\Longrightarrow)$:  Every countable model has a superminimal elementary end  extension. After checking the proof 
of this (for example,  in \cite[Theorem~2.1.12]{ksbook}), we see that if  $\MM_1$ is coded, then a superminimal extension $\MM_2$ can be constructed so that $\MM_2 \prec \MM$. 

\smallskip  $(\Longrightarrow)$: Let $\MM_2 \prec \MM$ be a superminimal extension 
of $\MM_1$. (It is not necessary that it be an end extension.) First note that $\MM_2$ is finitely generated, so it is coded. 
Let $b \in M$ code $M_2$; that is, $M_2 =  \{(b)_n : n < \omega\}$. 
Let $t_0(x), t_1(x),$ $ t_2(x), \ldots$ be a recursive list of all Skolem terms. 
Let $I \subseteq \omega$ be such that 
$i \in I$ iff $t_i(b) \in M_1$. Thus, $i \in I$ iff there is no Skolem term $t(y)$ 
such that $t(t_i(b)) = b$. Clearly, $I$ is recursive in $\tp(b)'$. 
Let $i_n$ be the $n$-th member of $I$. Since $\MM$ is arithmetically saturated,  
the set of formulas
$\{(x)_n = t_{i_n}(b) : n < \omega\}$ is realized in $\MM$ by, say, $a \in M$. 
Then $a$ codes $M_1$.  \qed

\bigskip

  {\sc Corollary 4.8}:   {\em Suppose that $\MM, \NN$ are countable, arithmetically saturated models 
and that $\alpha : \aut(\MM) \into \aut(\NN)$ is an isomorphism. 
If $\MM_0 \prec \MM$ and $\MM_0$ is coded, then $\wt{\alpha}(\MM_0)$ is coded in $\NN$.}

\bigskip

 {\it Proof}. Suppose that $\MM_0 \prec \MM$ is coded. By Lemma~4.7, let $\MM_1\prec \MM$ be a superminimal end extension of $\MM_0$. 
Then, by Lemma~4.3(2),   $\wt{\alpha}(\MM_1)$ is a superminimal end extension of $\wt{\alpha}(\MM_0)$,
so, by Lemma~4.7,    $\wt{\alpha}(\MM_0)$ is coded in $\NN$.   \qed

\bigskip

{\sc Lemma 4.9}: {\em Suppose that $\MM, \NN$ are countable,  arithmetically saturated models 
and that $\alpha : \aut(\MM) \into \aut(\NN)$ is an isomorphism. 
If $\MM_1, \MM_2 \preccurlyeq \MM$ are coded and $\ssy(\MM_1) \subseteq \ssy(\MM_2)$, then $\ssy(\wt{\alpha}(\MM_1)) \subseteq \ssy(\wt{\alpha}(\MM_2))$.}

\bigskip

{\it Proof}. This lemma  was proved in \cite[Coro.\@ 5.3]{ns} under the additional 
hypotheses that  $\MM$, $\NN$ are 2-Ramsey and that $\MM_1$, $\MM_2$ are finitely generated. The requirement that $\MM$, $\NN$ are $2$-Ramsey was needed only because it  also appeared in 
Proposition~4.2 of \cite{ns}, which is our Lemma~4.4 except that Proposition~4.2 has  the added hypothesis that 
$\MM$ and $\NN$ are $2$-Ramsey. 

Now suppose that $\MM_1, \MM_2$ are coded (as in  the lemma) and that $\ssy(\MM_1) \subseteq \ssy(\MM_2)$. Then, by Lemma~4.7, 
there are superminimal end extensions $\MM_3 \succ_{\sf end} \MM_1$ and 
$\MM_4 \succ_{\sf end} \MM_2$ and, therefore,   $\ssy(\MM_3) = \ssy(\MM_1) \subseteq \ssy(\MM_2) = \ssy(\MM_4)$.
Then,  
$\wt{\alpha}(\MM_3)$, $\wt{\alpha}(\MM_4)$ are finitely generated and, by  Lemma~4.3(2), 
are end extensions of  $\wt{\alpha}(\MM_1)$, $\wt{\alpha}(\MM_2)$, repectively. Thus, we that 
$\ssy(\wt{\alpha}(\MM_1)) = \ssy(\wt{\alpha}(\MM_3)) \subseteq \ssy(\wt{\alpha}(\MM_4)) = \ssy(\wt{\alpha}(\MM_2))$. \qed

\bigskip
The next lemma is a variation of Corollary~2.7 that takes place inside 
an arithmetically saturated model.

\bigskip

{\sc Lemma 4.10}: {\em Suppose that $\MM$ is an arithmetically saturated model, 
 $\MM_0 \prec \MM$ is coded and $X \in \ssy(\MM) \backslash \ssy(\MM_0)$. 
 Then there is a coded $\MM_1 \succ_{\sf cf} \MM_0$ such that $X \not\in \ssy(\MM_1)$ 
 and $\MM_1$ has the $\omega$-property.}  
 
 \bigskip
 
 {\it Proof.} The main concern is that the appropriate variant of Lemma~2.6 can be 
 proved. It is left to  the reader to show that the proof of Lemma~2.6 can be appropriately 
 modified. \qed
 
 \bigskip
 
 {\it Remark}: In this lemma, it is required that $\MM_1$ have the $\omega$-property.
 In fact, it will have the stronger property that there are a coded $\MM_2 \succ_{\sf end} 
 \MM_1$ and $a \in M_2$ such that $M_1 = \sup^{\MM_2}\{(a)_i : i < \omega\}$. 
 
 \bigskip

The next lemma is a variant of Corollary~2.10 that takes  place inside an arithmetically saturated 
model.

\bigskip

{\sc Lemma 4.11}: {\em Suppose that $\MM$ is an arithmetically saturated model, 
 $\MM_0 \prec \MM$ is coded and ${\mathfrak X} = \ssy(\MM_0)$.  
The following are equivalent$:$

$(1)$ $\Th(\MM) \not\in {\mathfrak X};$

$(2)$  for every tall coded $\MM_1 \prec \MM$ such that $\ssy(\MM_1) = {\mathfrak X}$,   
there is a coded $\MM_2\succcurlyeq_{\sf cf} \MM_1$ that is not generated by a set of indiscernibles.}

\bigskip

{\it Proof}. Let $T = \Th(\MM)$.

$(1) \Longrightarrow (2)$: Suppose that $T \not\in {\mathfrak X}$.  
Let $\MM_1 \prec \MM$ be tall and coded such that $\ssy(\MM_1) \subseteq {\mathfrak X}$. 
By Lemma~4.10, let $\MM_2 \prec \MM$ be such that 
$\MM_1 \prec_{\sf cf} \MM_2$, $T \not\in \ssy(\MM_2)$ and $\MM_2$ has 
the $\omega$-property. By Lemma~2.2, $\MM_2$ is uniformly $\omega$-logty. 
Since $T \not\in \ssy(\MM_2)$, $\MM_2$ is not recursively saturated, so, by Lemma~2.3, $\MM_2$ is not generated by a set of indiscernibles. 

\smallskip

$(2) \Longrightarrow (1)$: Suppose that $T \in {\mathfrak X}$.  
Let $\MM_1 \prec \MM$ be coded, recursively saturated and be such that 
$\ssy(\MM_1) \subseteq \ssy(\MM_0)$. Then,  $\MM_1$ is tall, and every countable $\MM_2 \succcurlyeq_{\sf cf}\MM_1$ is recursively saturated and, hence,  generated by a set of indiscernbles. \qed

\bigskip

The next corollary says  that  those elementary submodels that code their own theories are recognizable.

\bigskip

{\sc Corollary 4.12}: {\em  Suppose that $\MM, \NN$ are countable, arithmetically saturated models 
and that $\alpha : \aut(\MM) \into \aut(\NN)$ is an isomorphism. 
If $\MM_1 \preccurlyeq \MM$ and $\Th(\MM) \in \ssy(\MM_1)$, then 
$\Th(\NN) \in \ssy({\wt \alpha}(\MM_1))$.}

\bigskip

{\it Proof}. Suppose that $\MM_1 \preccurlyeq \MM$ and $\Th(\MM) \in \ssy(\MM_1)$. 
Let $a \in M$ and $\MM_0 = \scl(a)$ be such that  $\Th(\MM) \in \ssy(\MM_0)$. 
Since ${\wt \alpha}(\MM_0) \preccurlyeq {\wt \alpha}(\MM_1)$, 
it suffices to show that $\Th(\NN) \in \ssy({\wt \alpha}(\MM_0))$.

By $(2) \Longrightarrow (1)$ of Lemma~4.11, there is a tall coded $\MM_2$ such that $\ssy(\MM_2) = \ssy(\MM_0)$ and every coded $\MM_3 \preccurlyeq_{\sf cf} \MM_2$ 
is generated by a set of indiscernibles.

For a contradiction, suppose that $\Th(\NN) \not\in \ssy(\NN_0)$. 
Since $\NN_2$ is tall and coded and $\ssy(\NN_2) = \ssy(\NN_0)$, then by 
By $(1) \Longrightarrow (2)$ of Lemma~4.11, there is a tall coded 
$\NN_3 \succcurlyeq \NN_2$ that is generated by a set of indiscernibles. 
But then $\MM_3 \succcurlyeq \MM_2$ is tall, coded and  generated by a set of indiscernibles, which is a contradiction.
 \qed
 
 \bigskip

Corollary~4.12  has a $*$-version, so it can be generalized from models $\MM$ to 
models $(\MM,a)$ as in part (a) of the next lemma. Let us say that a model $\MM$ is a 
{\it self-coder} if $\tp(a) \in \ssy(\MM)$ whenever $a \in M$. Every recursively saturated model is a self-coder. Moreover, if 
 $\MM_0$ is countable, then, $\MM_0$ is a self-coder 
iff there is a countable, recursively saturated $\MM_1 \succ \MM_0$ such that 
$\ssy(\MM_1) = \ssy(\MM_0)$.
Part (b) of the next lemma says that 
self-coders are recognizable.
\bigskip

{\sc Corollary 4.13}: {\em  Suppose that $\MM, \NN$ are countable, arithmetically saturated models, 
$\alpha : \aut(\MM) \into \aut(\NN)$ is an isomorphism, and 
 $\MM_1 \preccurlyeq \MM$.

\begin{itemize}

\item[(a)] 
If  $a \in M_1$, $\wt{\alpha}(\scl(a)) = \scl(b)$ and $\tp(a) \in \ssy(\MM_1)$,  then $\tp(b) \in \ssy(\wt{\alpha}(\MM_1))$.

\item[(b)]  If $\MM_1$ is a self-coder, then $\wt{\alpha}(\MM_1)$ is a self-coder. \qed

\end{itemize}}

\bigskip

{\sc Lemma 4.14}: {\em Suppose that $\MM_0$ is countable. Then $\MM_0$ is recursively saturated iff 
each of the following holds$:$

\begin{itemize}

\item[(1)] $\MM_0$ is homogeneous$;$

\item[(2)] $\MM_0$ is a self-coder$;$

\item [(3)] if  $p(x)  \in \ssy(\MM_0)$ is a complete $1$-type for $\Th(\MM_0)$, then $p(x)$ is realized in 
$\MM_0$.

\end{itemize}}

\bigskip

{\it Proof}. $(\Longrightarrow)$: By Proposition~1.3, every countable, recursively saturated model is homogeneous and, as previously mentioned, every recursively saturated model is a self-coder. Lastly, observe that (3) is a straightforward 
consequence of recursive saturation. 

\smallskip

$(\Longleftarrow)$: For the converse, suppose that $\MM_0$ satisfies (1) -- (3). 
It follows from (2) that $\Th(\MM_0) \in \ssy(\MM_0)$. Thus, by Theorem~1.1,  there is 
a countable, recursively saturated $\MM_1 \equiv \MM_0$ such that $\ssy(\MM_1) = \ssy(\MM_0)$.  If $p(x)$ is a complete type, then (2) and (3) imply that $p(x)$ is realized 
in $\MM_0$ iff $p(x) \in \ssy(\MM_0)$. 
Thus, $\MM_0$ and $\MM_1$ realize exactly the same types. Since both are 
countable and homogeneous, 
they are isomorphic, so $\MM_0$ is also recursively saturated. \qed

\bigskip

The next lemma asserts that recursively saturated elementary submodels are recognizable.

\bigskip

{\sc  Lemma 4.15}: {\em Suppose that $\MM, \NN $ are countable, arithmetically saturated models, $\alpha : \aut(\MM) \into \aut(\NN)$ is an isomorphism, and $\MM_0 \preccurlyeq \MM$. If $\MM_0$ is recursively saturated, then $\wt{\alpha}(\MM_0)$ 
is recursively saturated.} 

\bigskip

{\it Proof}. 
Suppose that $\MM_0$ is recursively saturated, so that $\MM_0$ satisfies (1) -- (3) 
of Lemma~4.14. Let $\NN_0 = \wt{\alpha}(\MM_0)$. Then, $\NN_0$ is homogeneous by 
Lemma~4.5(1) and is a self-coder by Corollary~4.13(b). To get that $\NN_0$ is recursively saturated, it remains to show that 
it has the property described in (3) of Lemma~4.14. 
So, consider a complete type 
$q(x) \in \ssy(\NN_0)$, and let $b \in N$ realize $q(x)$ in $\NN$. 
Let $a \in M$ be such that $\wt{\alpha}(\scl(a)) = \scl(b)$. 
By Corollary~4.13(a), $\tp(a) \in \ssy(\MM_0)$, and, therefore, by Lemma~4.14(3), 
$\MM_0$ realizes $\tp(a)$. Let $a' \in M_0$ be such that $\tp(a') = \tp(a)$. 
Then, $\scl(a) \cong \scl(a')$, so by   Lemma~4.2, $\scl(b) \cong \wt{\alpha}(\scl(a'))$. 
Let $b' \in \wt{\alpha}(\scl(a'))$ be such that $\tp(b') = \tp(b)$. But then 
$b'  \in N_0$ and $\tp(b') = q(x)$.  \qed

\bigskip

We can now complete the proof of Theorem~4. Suppose that $\MM, \NN$ are countable,
arithmetically saturated models and that $\alpha : \aut(\MM) \into \aut(\NN)$ is an 
isomorphism. Let $T = \Th(\MM)$ (so that $T'$ is the Turing-jump of $T$).  It suffice to prove that $T' \leq_T \Th(\NN)'$.

Since $T \in \ssy(\MM_0)$ for every recursively saturated $\MM_0 \preccurlyeq \MM$, 
we get from  $(2) \Longrightarrow (1)$ of Lemma~3.1, that every recursively saturated 
$\MM_0 \preccurlyeq \MM$ has an elementary substructure   $\MM_1 \prec \MM_0$ such that 
$\Lt(\MM_1) \cong {\mathcal D}'(T')$. Then Lemma~4.15 implies that 

\begin{itemize}
\item[$(*)$]
  every recursively saturated $\NN_0 \preccurlyeq \NN$ has an elementary \\ substructure 
 $\NN_1 \prec \NN_0$ such that 
$\Lt(\NN_1) \cong {\mathcal D}'(T')$. 

\end{itemize}
\noindent
Now suppose, for a contradiction, that $T' \not\leq_T \Th(\NN)'$.
There is a countable Scott set ${\mathfrak X} \subseteq \ssy(\NN)$ such that $\Th(\NN) \in {\mathfrak X}$ and 
$X' \leq_T \Th(\NN)'$ for every $X \in {\mathfrak X}$. (See \cite[Theorem~VIII.2.17]{simp}.) 
Hence, $T' \not\leq_T X'$ for every $X \in {\mathfrak X}$.  By Theorem~1.1, let $\NN_0$ be recursively saturated model of $\Th(\NN)$  such that 
$\ssy(\NN_0) = {\mathfrak X}$. Since ${\mathfrak X} \subseteq \ssy(\NN)$, we can 
assume that $\NN_0 \prec \NN$. Then, by $(1) \Longrightarrow (2)$ of Theorem~3.1,
there is no $\NN_1 \prec \NN_0$ such that $\Lt(\NN_1) \cong {\mathcal D}'(T')$, 
contradicting $(*)$ and  completing the proof of Theorem~4. \qed

\bigskip

We next prove Corollary~5. We will actually prove something somewhat stonger since 
we will obtain completions $T_0,T_1,T_2, \ldots$ such that each $\MM_i$ that we obtain 
is a model of $T_i$.

Let ${\mathfrak X}_0$ be a Scott set 
that is enumerated by some arithmetical $X$. Apply Theorem~1.2 to get, for each $n < \omega$, 
a completion $T_n$  of \pa\  such that $T_n \equiv_n X^{(2n)}$ 
and $\rep(T_n) = {\mathfrak X}_0$. For any countable jump ideal ${\mathfrak X}$, Theorem~1.1 implies that there 
are models $\MM_n \models  T_n$ such that $\ssy(\MM_n) = {\mathfrak X}$. 
Each $\MM_n$ is arithmetically saturated since ${\mathfrak X}$ is a jump ideal. 
By Theorem~4, $\aut(\MM_i) \not \cong \aut(\MM_j)$ whenever $i < j < \omega$.

\bigskip


\bigskip

{\bf \S5.\@ The Proof of Theorem 6.} Lemma~4.15 showed that recursively saturated elementary submodels are recognizable. It is also the case that arithmetically saturated elementary submodels are 
recognizable. 

\bigskip

{\sc Lemma 5.1}:  {\em Suppose that $\MM, \NN $ are countable, arithmetically saturated models, $\alpha : \aut(\MM) \into \aut(\NN)$ is an isomorphism, and $\MM_0 \preccurlyeq \MM$. If $\MM_0$ is arithmetically  saturated, then $\wt{\alpha}(\MM_0)$ 
is arithmetically saturated.} 

\bigskip

{\it Proof}. There are several ways to prove this. One way is to use \cite[Coro.\@ 5.4]{kkk}
that asserts:  If $\MM_0$ is recursively saturated, then $\MM_0$ is arithmetically saturated iff there is $g \in \aut(\MM_0)$ and an open $H < \aut(\MM_0)$ such that for all 
$f \in \aut(\MM_0)$, 
$f^{-1}gf \not\in H$. We leave it to the reader to complete this proof. \qed 

\bigskip

 Suppose that $T$ is any completion of \pa. Theorem~1.1 implies that there is a (necessarily unique) recursively saturated model 
 $\MM \models T$ such that $\ssy(\MM)$ is the smallest jump ideal to which $T$ belongs. 
This model is arithmetically saturated and is elementarily embeddable in every 
arithmetically saturated model of $T$. We refer to this $\MM$ as the {\it minimal 
arithmetically saturated} model of~$T$. 
Minimal 
arithmetically saturated models are recognizable.

\bigskip

{\sc Corollary 5.2}: {\em Suppose that $\MM, \NN $ are countable, arithmetically saturated models, $\alpha : \aut(\MM) \into \aut(\NN)$ is an isomorphism, and $\MM_0 \preccurlyeq \MM$. If $\MM_0$ is a minimal arithmetically  saturated model, then so is $\wt{\alpha}(\MM_0)$.}  \qed

\bigskip

Recall the following key result of \cite[Theorem~3.8]{n}.

\bigskip 

{\sc  Lemma 5.3}: {\em Suppose that $\MM, \NN $ are countable, arithmetically saturated models, $\alpha : \aut(\MM) \into \aut(\NN)$ is an isomorphism, and $1 \leq n < \omega$. If $(\mathbb N,\rep(\Th(\MM))) \models {\sf RT^n_2}$, then  $(\mathbb N,\rep(\Th(\NN))) \models {\sf RT^n_2}$.} 

\bigskip

The following is a rephrasing of Lemma~5.3 when $n=3$.

\bigskip

{\sc Corollary 5.4}: {\em Suppose that $\MM, \NN $ are countable, arithmetically saturated models, $\alpha : \aut(\MM) \into \aut(\NN)$ is an isomorphism. If $\MM_0 \prec \MM$ and $\NN_0 \prec \NN$ are the prime elementary submodels and $\ssy(\MM_0)$ is a jump ideal, 
then $\ssy(\NN_0)$ is a jump ideal.} \qed

\bigskip

The $*$-version of the previous corollary implies that finitely generated models whose 
standard systems are jump ideals are recognizable. But this is also so 
for coded models.

\bigskip

 {\sc  Lemma 5.5}: {\em Suppose that $\MM, \NN $ are countable, arithmetically saturated models, $\alpha : \aut(\MM) \into \aut(\NN)$ is an isomorphism, and $\MM_0 \prec \MM$. If $\MM_0$ is coded and $\ssy(\MM_0)$ is a jump ideal, 
then $\ssy(\wt{\alpha}(\MM_0))$ is a jump ideal.}

\bigskip

{\it Proof}. By Lemma~4.7, let $\MM_1\prec \MM$ be a superminimal end extension of 
$\MM_0$. Then $\MM_1$ is finitely generated and $\ssy(\MM_1) = \ssy(\MM_0)$. 
Also, $\wt{\alpha}(\MM_1)$ is finitely generated, $\wt{\alpha}(\MM_1)
\succ_{\sf end} \wt{\alpha}(\MM_0)$ (by Lemma~4.3(2)) and $\wt{\alpha}(\MM_1)$ is finitely generated. Therefore, $\ssy(\wt{\alpha}(\MM_1))$ is a jump ideal and so is 
$\ssy(\wt{\alpha}(\MM_0))$. \qed

\bigskip

Given an arithmetically saturated model $\MM$, we let $j(\MM)$ be the cardinality of the set of all jump ideals ${\mathfrak X}$ for which there are 
$\MM_0 \preccurlyeq \MM_1 \preccurlyeq \MM$ such that 
$\MM_0$ is coded, $\ssy(\MM_0) = {\mathfrak X}$ and $\MM_1$ is a minimal arithmetically 
saturated model. 

\bigskip

 {\sc  Lemma 5.6}: {\em Suppose that $\MM, \NN $ are countable, arithmetically saturated models and $\aut(\MM) \cong \aut(\NN)$. Then, $j(\MM) = j(\NN)$.} 
 
 \bigskip
 
 {\it Proof}. Immediate from Lemmas~4.9, 5.2 and  5.5. \qed
 
 \bigskip
 
 We will say that a completion $T$ of \pa\ is {\it tight} if there is $m < \omega$ such that whenever $\MM$ is a countable, arithmetically saturated model of $T$, then $j(\MM) = m$.
 If $T$ is tight, then we let $j(T)$ be that $m$. 
 
 \bigskip

We now prove Theorem~6. Fix $n < \omega$. Our goal is to obtain 
 recursively equivalent completions $T_0,T_1, \ldots, T_n$ such that for each $i \leq n$, 
 $T_i$ is tight and $j(T_i) = n+1-i$. Clearly, by Lemma~5.6, these theories 
will suffice to confirm Theorem~6.

Consider any $A_0 \subseteq \omega$ and then let 
$A_0, A_1,A_2, \ldots , A_n \subseteq \omega$ be such that  
$A_0 <_a A_1 <_a A_2 <_a \cdots  <_a A_n$ and whenever $A$ is such that $A_0\leq _a A \leq_a A_n$, then there is a unique $i \leq n$ such that $A \equiv_a A_i$. 
That there are $A_1,A_2, \ldots, A_n$ follows from a very special case of the theorem of Harding \cite{hard}.
For each $i \leq n$, let ${\mathcal B}_i$ be the jump ideal $\{B \subseteq \omega : B \leq_a A_i\}$, and then let $X_i \in {\mathcal B}_i$ and ${\mathfrak X}_i \subseteq {\mathcal B}_i$ be such that ${\mathfrak X}_i$ is a Scott set, $A_i \in {\mathfrak X}_i$ and $X_i$ enumerates
${\mathfrak X}_i$.  (See \cite[Theorem~VIII.2.17]{simp}.) 
 Let 
$X = X_0 \oplus X_1 \oplus \cdots \oplus X_n \oplus A_n$, so ${\mathcal B}_n$ is the smallest jump ideal containing $X$. 
By Theorem~1.2, let $T_i$ be a completion 
of $\pa$ such that $T_i \equiv_T X$ and $\rep(T_i) = {\mathfrak X}_i$.  
Each $T_i \in {\mathcal B}_n$, so by Theorem~1.1, there is a countable, recursively saturated $\MM_i \models T_i$ such that $\ssy(\MM_i) = {\mathcal B}_n$. 
Clearly, $\MM_i$ is a minimal arithmetically saturated model of $T_i$. 
It is clear that $j(\MM_i) \leq n+1-i$ as that is the cardinality of the set 
$\{{\mathcal B}_i, {\mathcal B}_{i+1},\ldots,{\mathcal B}_n\}$, which is the set of all jump ideals ${\mathfrak X}$ such that  $\rep(T_i) \subseteq {\mathfrak X} \subseteq \ssy(\MM_i)$. Moreover, if $\MM$ is any arithmetically saturated model of $T_i$, then $j(\MM) \leq n+1-i$. 

Thus, it remains to show 
that 
$j(T_i) \leq n+1-i$; that is, we must show that whenever $i \leq j \leq n$, then there 
is a coded elementary submodel of $\MM_i$ whose standard system is ${\mathcal B}_j$.
This is a consequence of the next lemma, which we state separately since it has its own interest.
 \qed

\bigskip
 
 {\sc Lemma 5.7:} {\em Suppose that $\MM$ is arithmetically saturated,  
 ${\mathfrak X}_0$ is a Scott set enumerated by some $X \in \ssy(\MM)$ and $\rep(\Th(\MM)) \subseteq {\mathfrak X}_0$.
 Then there is a finitely generated $\MM_0 \prec \MM$ such that 
 $\ssy(\MM_0) = {\mathfrak X}_0$.}

\bigskip

{\it Proof}.  For $n < \omega$, let ${\mathcal L}_n = {\mathcal L}_\pa \cup \{c_0,c_1, \ldots, c_{n-1}\}$, where the $c_i$'s are new and distinct constant symbols.
Thus, ${\mathcal L}_0 = {\mathcal L}_\pa$. Let ${\mathcal L} = {\mathcal L}_\pa \cup \{c_0,c_1, c_2,\ldots\}$. 

Suppose that $X \subseteq \omega$ enumerates ${\mathfrak X}_0$. (For the time being, no other assumptions are being made.) To simlify notation, let $X_n = (X)_n$.

We will construct a complete Henkin ${\mathcal L}$-theory $T \supseteq \Th(\MM)$,
and then let $\MM_0$ be the Henkin model of $T$. The theory $T$ will be the union 
of an increasing chain $T_0 \subseteq T_1 \subseteq T_2 \subseteq \cdots$ of theories such that 
for each $n < \omega$,
\begin{itemize}

\item[(1)] $T_n$ is an ${\mathcal L}_n$-theory and $T_n \subseteq \Pi_n \cup \Sigma_n$;

\item[(2)] for every ${\mathcal L}_n$-sentence $\sigma \in \Pi_n$, either $\sigma \in T_n$
or $\neg\sigma \in T_n$;

\item[(3)] $\{(b_{2n})_i = 0 : i \in X_n\} \cup \{(b_{2n})_i = 1 : i \in \omega \backslash X_n\} \subseteq T_{2n+1}$;

\item[(4)] the sentence $\exists x \varphi_n(x) \rightarrow \varphi_n(b_{2n+1})$ is in $T_{2n+2}$.

\item[(5)] $T_n \in {\mathfrak X}_0$;

\item[(6)] $T_n \cup \Th(\MM)$ is consistent.

\end{itemize}

Suppose that we have obtained the sequence $T_0 \subseteq T_1 \subseteq T_2 \subseteq \cdots$ satisfying (1) -- (5) and that $T$ is its union. 

It follows from (1) and (6)  that $T$ is an ${\mathcal L}$-theory, 
from (2) that $T$ is complete, 
from (2) and (6)  that $T \supseteq \Th(\MM)$, from (4) that $T$ is a Henkin theory,
from (3) that $\rep(T) \supseteq {\mathfrak X}_0$, and from (5) that $\rep(T) \subseteq {\mathfrak X}_0$. Thus, we can let $\MM_1$ be  the ${\mathcal L}_\pa$-reduct of  the Henkin model of $\Th(\MM)$. Hence,  $\MM_1 \equiv \MM$ 
and $\ssy(\MM_1) = {\mathfrak X}_0$. 

The construction of the sequence of theories proceeds by recursion.

\smallskip

Let $T_0 = \Pi_0 \cap \Th(\MM)$.

\smallskip

Suppose that $n = 2k+1$. Let $T'_n = T_{n-1} \cup \{(b_{2k})_i = 0 : i \in X_k\} \cup \{(b_{2k})_i = 1 : i \in \omega \backslash X_k\}$. Clearly, $T'_n \in {\mathfrak X}_0$ and 
$T'_n \cup \Th(\MM)$ is consistent. Let  $T''_n \cup ((\Pi_n \cup \Sigma_n) \cap \Th(\MM)) $.
Then,  $T''_n \in {\mathfrak X}$ and 
$T''_n$ is consistent, so there is a complete ${\mathcal L}_n$-theory $T'''_n \supseteq T''_n$.
Let $T_n = T'''_n \cap (\Pi_n \cup \Sigma_n)$.

\smallskip

Suppose that $n = 2k+2$. Let $T'_n = T_{n-1} \cup \{\exists x \varphi_n(x) \rightarrow \varphi_k(b_{2k+1})\}$. Clearly, $T'_n \in {\mathfrak X}$ and 
$T'_n \cup \Th(\MM)$ is consistent. Now obtain $T_n$ from $T'_n$ exactly as was done in the previous case. 

\smallskip

We have shown how to obtain a model $\MM_1 \equiv \MM$ such that 
$\ssy(\MM_1) = {\mathfrak X}_0$. (Thus, we have proved one of Scott's theorems.) 
We next indicate how the construction can be refined to a more effective one that 
yields a coded $\MM_1 \prec \MM$.

First, we recall a theorem of Marker \cite{mm} that says that if a Scott set has an enumeration, 
then it has an {\it effective} enumeration. (More details about this, including a definition,
can be found in \cite[Chap.\@ 19]{akbook} or in \cite{knight}.) It can be checked that 
if a Scott set has an enumeration in a jump ideal ${\mathfrak X}$, then it also has an effective enumeration in ${\mathfrak X}$. Thus, we can assume that ${\mathfrak X}_0$ 
is effectively enumerated by $X \in \ssy(\MM)$. 

It can be checked that the construction of the sequence of theories 
is recursive in $X \oplus T$. Let $\Phi(x)$ be the set of $1$-ary ${\mathcal L}_\pa$-formulas 
that are obtained by replacing each occurrence of $c_n$ by $(x)_n$. 
Then, $\Phi(x) \in \ssy(\MM)$ and is consistent with $\Th(\MM)$, so there is some 
$a \in M$ that satisfies $\Phi(x)$. Thus, without loss of generality, $a$ codes $\MM_1$.
By Lemma~4.7, let $\MM_0$ be a superminimal end extension of $\MM_1$ such that 
$\MM_0 \prec \MM$. \qed

\bigskip

We end this section by asking if Theorem~6 can be improved. 

\bigskip

{\it Question} 5.8:  Are there  infinitely many, recursively equivalent completions $T_0,T_1, T_2, \ldots$ of $\pa$ such that whenever 
$i < j < \omega$ and $\MM_i,\MM_j$ are countable, arithmetically saturated models of 
$T_i,T_j$, respectively, then  $\aut(\MM_i) \not\cong \aut(\MM_j)$?

\bigskip

In the absence of a positive answer to the previous question (or even following an unlikely 
negative answer), we could still ask the next question. 

\bigskip

{\it Question} 5.9:  Are there infinitely many countable, arithmetically saturated models 
no two of which have isomorphic automorphism groups but
 all of which have the same standard systems and 
recursively equivalent theories?

\bigskip


{\bf \S6.\@ Some Additional Results.}   A consequence of Theorems~1 and~4  is that there is a set ${\mathcal T}$ of $2^{\aleph_0}$ completions  of \pa\ such that if 
$\MM$, $\NN$ are  nonisomorphic countable arithmetically saturated
models and $\Th(\MM), \Th(\NN) \in {\mathcal T}$, then 
$\aut(\MM) \not\cong \aut(\NN)$. The next corollary is a strengthening of this.

\bigskip

{\sc Corollary 6.1}: {\em Suppose that ${\mathfrak X}$ is a countable Scott set. 
There is a set ${\mathcal T}$ of $2^{\aleph_0}$ completions  of \pa\ such that 
every $T \in \mathcal T$ has $\rep(T) = {\mathfrak X}$ and if 
$\MM$, $\NN$ are  nonisomorphic countable arithmetically saturated
models such that $\Th(\MM), \Th(\NN) \in {\mathcal T}$, then 
$\aut(\MM) \not\cong \aut(\NN)$.} 

\bigskip

{\it Proof.} Let $X \subseteq \omega$ enumerate ${\mathfrak X}$. For every $Y \geq_T X$, 
apply Theorem~1.1 to get a completion  $T_Y \equiv_T Y$ such that $\rep(T_Y) = {\mathfrak X}$. Then ${\mathcal T} = \{T_Y : Y \geq_T X\}$ has cardinality $2^{\aleph_0}$. \qed

\bigskip

The proof of the main result of \cite{ns} showed that if  $\MM, \NN$ are  saturated 
models and 
  $\aut(\MM) \cong \aut(\NN)$, then $ \Th(\MM)^{(\omega)} \equiv_a \Th(\NN)^{(\omega)}$. The conclusion can be improved in the manner of Theorem~4 with essentially the same proof as Theorem~4.

\bigskip

{\sc Theorem 6.2}:  {\em If $\MM, \NN$ are  saturated 
models and 
  $\aut(\MM) \cong \aut(\NN)$, then $ \Th(\MM)' \equiv_T \Th(\NN)'$.} 
  
  \bigskip
  
  Theorem~6 also has a version for saturated models. 
  
  \bigskip

{\sc Theorem 6.3}: { \em For each $n < \omega$, there are  recursively equivalent 
completions $T_0,T_1, \ldots, T_n$ of\  \pa\ such that 
whenever $i < j \leq n$ and $\MM_i \models T_i$, $\MM_j \models T_j$ are saturated models, then 
$\aut(\MM_i) \not\cong \aut(\MM_j)$.} 

\bigskip

{\it Proof}. Use the same theories as in the proof of Theorem~6. \qed

\bigskip

There is another approach to proving the previous two theorems using the corresponding results for arithmetically saturated models and the following 
lemma, the proof of which we omit but which the reader should be able to work out.

\bigskip

{\sc Lemma 6.4}: {\em Suppose that  $\MM, \NN$ are  saturated 
models and that 
  $\aut(\MM) \cong \aut(\NN)$. If $\MM_0 \prec \MM$ and $\NN_0 \prec \NN$ are minimal arithmetically saturated models, then $\aut(\MM_0) \cong \aut(\NN_0)$.}

\bibliographystyle{plain}

\end{document}